# PYTHAGOREAN FUZZY GRAPHS: SOME RESULTS

**Rajkumar Verma [1, 2], José M. Merigó[1], Manoj Sahni[3]**

[1]*Department of Management Control and Information Systems, University of Chile,*
*Av. Diagonal Paraguay 257, Santiago-8330015, CHILE*

[2]*Department of Applied Sciences, Delhi Technical Campus,*
*(Affiliated to Guru Govind Singh Indraprastha University, Delhi)*
*28/1, Knowledge Park-III, Greater Noida-201306*
*Uttar Pradesh, INDIA*

[3]*Department of Mathematics, School of Technology,*
*Pandit Deendayal Petroleum University,*
*Gandhinagar-382007*
*Gujarat, INDIA*
*rverma@fen.uchile.cl, jmerigo@fen.uchile.cl, manoj_sahani117@rediffmail.com*

## ABSTRACT

Graph theory has successfully used to solve a wide range of problems encountered in diverse fields such as medical sciences, neural networks, control theory, transportation, clustering analysis, expert systems, image capturing, and network security. In past few years, a number of generalizations of graph theoretical concepts have developed to model the impreciseness and uncertainties in graphical network problems.

A Pythagorean fuzzy set is a powerful tool for describing the vague concepts more precisely. The Pythagorean fuzzy set-based models provide more flexibility in handling the human judgment information as compared to other fuzzy models. The objective of this paper is to apply the concept of Pythagorean fuzzy sets to graph theory. This work introduces the notion of Pythagorean fuzzy graphs (PFGs) and describes a number of methods for their construction. We then define some basic operations on PFGs and prove some of their important properties. The work also discusses the notion of isomorphism between Pythagorean fuzzy graphs with a numerical example. Further, we introduce the concept of the strong Pythagorean fuzzy graph and the complete Pythagorean fuzzy graph. In addition, the paper also proves some results on self-complementary, self-weak complementary with Pythagorean fuzzy strong graphs and Pythagorean fuzzy complete graphs.

**1. Introduction**

Zadeh [53] introduced the notion of fuzzy sets to model the uncertainty or vague concepts by assigning a membership degree corresponding to each element whose range is in between 0 and 1. Since the pioneering work of Zadeh, the fuzzy set theory has been used in different disciplines such as management sciences, engineering, mathematics, social sciences, statistics, signal processing, artificial intelligence, automata theory, medical and life sciences. In 1982, Atanassov [6] generalized the idea of fuzzy sets and introduced a new set theory called the 'intuitionistic fuzzy sets'. In the intuitionistic fuzzy set, each element has degrees of membership and non-membership whose sum lies between 0 and 1. In last three decades, the intuitionistic fuzzy set theory has been widely studied and a great number of applications have been developed in various fields including decision making [23, 42, 43], medical diagnosis [13, 40, 44], market prediction [18], clustering analysis [48, 49, 58] and pattern recognition [10-12, 22].

Graph theory is an important branch of applied mathematics and has numerous applications in different disciplines including computer science, economics, social sciences, chemistry, physics, system analysis, neural networks, electrical engineering, control theory, transportation, architecture and communication [14, 17]. However, in many realistic situations, some aspect of a graph-theoretic problem may be uncertain and cannot be represented by Euler's graph.

To handle such type of situations, in 1975, Rosenfeld [35] generalized the Euler's graph theory and proved the basic results on fuzzy graphs (FGs). Bhattacharya [7] made some comments on FGs and established some connectivity concepts regarding fuzzy cutnotes and fuzzy brides. Further, Bhutani [9] studied the automorphisms on fuzzy graphs and proved a number of properties connected with the complete fuzzy graph. Mordeson and Peng [28] defined some basic operations on fuzzy graphs. Mcallister [26] discussed theoretical and computational aspects of fuzzy intersection graphs by matrix representation. Bhattacharya and Suraweera [8] proposed an algorithm to find the connectivity of a pair of nodes in a fuzzy graph. In 1993, Mordeson [27] studied fuzzy line graph and proved the necessary and sufficient condition for a fuzzy graph to be isomorphic to its corresponding line graph. In 2001, Mordeson and Nair [29] defined the complement of a fuzzy graph. Later, Sunitha and Kumar [41] proposed a modified definition of the complement of a fuzzy graph and proved some properties on self-complementary fuzzy graphs. The cofuzzy graphs were studied by Akram [1] in 2011. Samanta and Pal [38] discussed the notion of fuzzy planar graphs and made a comparative study between Kuratowski's graphs and fuzzy planar graph.

In 1994, Shannon and Atanassov [39] further generalized the fuzzy graph theory and proposed the intuitionistic fuzzy graphs (IFGs) with some fundamental results. Karunambigai and Parvathi [20] analyzed the properties of minmax intuitionistic fuzzy graphs. Parvathi et al. [30] defined some operations on IFGs and studied their properties. Karunambigai et al. [21] discussed the constant-IFGs and totally constant-IFGs. In 2012, Akram and Davvaz [4] studied the strong-IFGs and proved some of their properties. Further, Akram and Dudek [2] proposed intuitionistic fuzzy hypergraphs and discussed their applications. Alshehri and Akram [5] developed the notion of multigraphs, planar graphs, and dual graphs under intuitionistic fuzzy environment. Sahoo and Pal [36] defined different types of product operations on IFGs. Recently, Sahoo and Pal [37] have proposed the idea of intuitionistic fuzzy tolerance graph, intuitionistic fuzzy $\phi$-tolerance graph and discussed their applications. Due to their flexibility and applicability, a number of researchers have been started work on intuitionistic fuzzy graph theory and proved a number of interesting results.

More recently, Yager [50, 52] has proposed the Pythagorean fuzzy set (PFS) as an effective tool for handling/modeling the uncertainty or vague information more adequately in real-world situations. In PFSs, the sum of squares of the degrees of membership and nonmembership is less than or equal to 1. For example, if a decision maker provides the membership degree 0.6 and nonmembership degree 0.7 in his evaluation, then this situation cannot handle by intuitionistic fuzzy set theory because of 0.6+0.7 > 1. However, it is easily observed that $0.6^2+0.7^2 < 1$, that is to say, the Pythagorean fuzzy set (PFS) is capable to represent this evaluation information. In other words, the PFSs are more powerful to handle problems in uncertain situations. Under PFS environment, many researchers have started work in different directions and obtained various significant results. Yager and Abbasov [51] established a relation between Pythagorean membership degrees and complex numbers and proved that the Pythagorean degrees are a subclass of complex numbers. Zhang and Xu [59] made a detailed study on Pythagorean fuzzy sets and proposed an extension of the TOPSIS method with Pythagorean fuzzy information. Peng and Yang [31] defined the division and subtraction operations on PFNs. Reformat and Yager [34] studied the collaborative-based recommender systems under Pythagorean fuzzy environment. Dick et al. [15] discussed some operations on Pythagorean and complex fuzzy sets. Ma and Xu [25] defined symmetric Pythagorean fuzzy weighted geometric/ averaging operators and studied their applications in multicriteria decision making. Peng and Yang [32] developed MABAC method under Pythagorean fuzzy environment. Zhang [57] discussed multicriteria Pythagorean fuzzy decision analysis using hierarchical QUALIFLEX approach with the closeness index based ranking methods. Zeng et al. [55] proposed a hybrid method for solving Pythagorean fuzzy multiple-criteria decision-making problems. Peng et al. [33] introduced some Pythagorean fuzzy information measures and discussed their applications. Further, Zeng [54] developed a new method to solve Pythagorean fuzzy multiattribute group decision making problems. In recent years, many researchers [24, 45-47, 56] have developed a number of aggregation operators to aggregate different PFNs.

Intuitionistic fuzzy graphs have been successfully applied in solving many problems connected with different areas [3, 19]. But there are many problems in real life which cannot be represented adequately by intuitionistic fuzzy graphs. So we need a more general graph theory to tackle such type of situations.

The aim of this work is to develop the graph-theoretic concepts under Pythagorean fuzzy environment. For doing so, we first propose the notion of Pythagorean fuzzy relation (PFR) and Pythagorean fuzzy graph (PFG) as a further generalization of FGs and IFGs. We then define some basic operations on PFGs and prove some their important properties. We also study the isomorphism and weak isomorphism between PFGs. Further, the work proposes the idea of the strong Pythagorean fuzzy graph and complete Pythagorean fuzzy graph. In addition, we also prove some results on self-complementary and self-weak complementary with Pythagorean fuzzy strong graphs and Pythagorean fuzzy complete graphs.

This paper is organized as follows: Section 2 describes some prerequisite material on graph theory, intuitionistic fuzzy graph theory, and Pythagorean fuzzy sets. Section 3 proposes the idea of Pythagorean fuzzy relation, Pythagorean fuzzy graph and some basic operations including

Cartesian product, composition, union, join and complement on Pythagorean fuzzy graphs. Some important properties of different operations on Pythagorean fuzzy graphs also proved in this section with illustrative numerical examples. In Section 4 we define an isomorphism between Pythagorean fuzzy graphs. Section 5 introduces the notion strong Pythagorean fuzzy graph and complete Pythagorean fuzzy graph. Further, a number of propositions are proved on strong Pythagorean fuzzy graphs and complete Pythagorean fuzzy graphs. Section 6 concludes the paper with some future directions.

## 2. Preliminaries

This section presents a brief review of graph theory, intuitionistic fuzzy graphs, and Pythagorean fuzzy sets, which will be used for further development.

### 2.1 *Some basic definitions in graph theory*

*Definition* **2.1.1** (*Graph*) [16]: A graph $G = (V, E)$ consists of two sets $V$ and $E$. The elements of $V$ are called vertices and the elements $E$ are called edges (formed by a pair of vertices). Two vertices $u$ and $v$ in $G$ are said to be adjacent if $(u, v)$ is an edge in $G$.

*Definition* **2.1.2** (*Simple Graph*) [16]: A graph is simple if it has no parallel edges or self-loops.

*Definition* **2.1.3** (*Complete Graph*) [16]: A complete graph is a simple graph in which every pair of distinct vertices is connected by an edge.

*Definition* **2.1.4** (*Complementary Graph*) [16]: The complement of a simple graph $G$ is the graph $\bar{G}$ that has the same vertices as $G$ such that any two vertices are adjacent in $\bar{G}$ if and only if they are not adjacent in $G$.

*Definition* **2.1.5** (*Cartesian product*) [16]: The Cartesian product of two simple graphs $G_1 = (V_1, E_1)$ and $G_2 = (V_2, E_2)$ is a graph defined by $G = G_1 \times G_2 = (V, E)$ with $V = V_1 \times V_2$ and

$$E = \{(u, u_2)(u, v_2) : u \in V_1, u_2 v_2 \in E_2\} \cup \{(u_1, v)(v_1, v) : v \in V_2, u_1 v_1 \in E_1\}.$$

*Definition* **2.1.6** (*Composition*) [16]: The composition of two simple graphs $G_1 = (V_1, E_1)$ and $G_2 = (V_2, E_2)$ is a graph defined by $G_1[G_2] = (V_1 \times V_2, E^\bullet)$ where $E^\bullet = E \cup \{(u_1, u_2)(v_1, v_2) \mid u_1 v_1 \in E_1, u_2 \neq v_2\}$ and $E$ is defined as in $G_1 \times G_2$. Note that $G_1[G_2] \neq G_2[G_1]$.

***Definition* 2.1.7 (*Union*)** [16]: The union of two simple graphs $G_1 = (V_1, E_1)$ and $G_2 = (V_2, E_2)$ is a simple graph defined by $G = G_1 \cup G_2 = (V_1 \cup V_2, E_1 \cup E_2)$.

***Definition* 2.1.8 (*Join*)** [16]: The join of two simple graphs $G_1 = (V_1, E_1)$ and $G_2 = (V_2, E_2)$ is denoted by $G_1 + G_2 = (V_1 \cup V_2, E_1 \cup E_2 \cup E')$, where $E'$ represents the set of all edges joining the nodes of $V_1$ and $V_2$ and assume that $V_1 \cap V_2 = \phi$.

***Definition* 2.1.9 (*Isomorphism*)** [16]: Two graphs $G_1 = (V_1, E_1)$ and $G_2 = (V_2, E_2)$ are isomorphic, written $G_1 \cong G_2$, if there are bijections $\theta : V(G_1) \to V(G_2)$ and $\psi : E(G_1) \to E(G_2)$ such that $\phi_{G_1}(s) = uv$ if and only if $\phi_{G_2}(\psi(s)) = \theta(u)\theta(v)$. Such a pair of mapping is called isomorphism between $G_1$ and $G_2$.

## 2.2 Intuitionistic fuzzy graphs: Basic results

Atanassov [6] extended the fuzzy set to the IFS and defined as follows:

***Definition* 2.2.1 (*Intuitionistic Fuzzy Set*)** [6]: An intuitionistic fuzzy set $\tilde{P}$ defined in a finite universe of discourse $U = (u_1, u_2, \ldots, u_n)$ is given by

$$\tilde{P} = \{\langle u, \varsigma_{\tilde{P}}(u), \xi_{\tilde{P}}(u)\rangle \mid x \in U\}, \tag{1}$$

where

$$\varsigma_{\tilde{P}}(u) : U \to [0,1], \xi_{\tilde{P}}(u) : U \to [0,1] \text{ and } 0 \leq \varsigma_{\tilde{P}}(u) + \xi_{\tilde{P}}(u) \leq 1 \quad \forall u \in U. \tag{2}$$

Here the numbers $\varsigma_{\tilde{P}}(u)$ and $\xi_{\tilde{P}}(u)$, respectively, denote the degree of membership and degree of non-membership of $u \in U$ in $\tilde{P}$.

For each intuitionistic fuzzy set $\tilde{P}$ in $U$, the intuitionistic fuzzy index (hesitation degree) can be defined as

$$\eta_{\tilde{P}}(u) = 1 - \varsigma_{\tilde{P}}(u) - \xi_{\tilde{P}}(u) \quad \forall u \in U. \tag{3}$$

If $\eta_{\tilde{P}}(u) = 0 \ \forall \ u \in U$, then IFS $\tilde{P}$ becomes fuzzy set. Therefore, fuzzy sets are special cases of IFSs.

Shannon and Atanassov [39] defined intuitionistic fuzzy graphs as an extension of fuzzy graph theory.

***Definition* 2.2.2 (*Intuitionistic Fuzzy Relation*)** [39]: Let $S$ and $T$ be two arbitrary finite non-empty sets. An intuitionistic fuzzy relation (IFR) $\tilde{R}(S \to T)$ in the universe $S \times T$ is an intuitionistic fuzzy set of the form:

$$\tilde{R} = \{\langle (s,t), \varsigma_{\tilde{R}}(s,t), \psi_{\tilde{R}}(s,t) \rangle \mid s \in S, t \in T\}, \tag{4}$$

where $\quad \varsigma_{\tilde{R}}(s,t): S \times T \to [0,1]$ and $\xi_{\tilde{R}}(s,t):(s,t): S \times T \to [0,1], \tag{5}$

and $\quad 0 \leq \varsigma_{\tilde{R}}(s,t) + \xi_{\tilde{R}}(s,t) \leq 1 \quad \forall (s,t) \in S \times T. \tag{6}$

Here the numbers $\varsigma_{\tilde{R}}(s,t)$ and $\xi_{\tilde{R}}(s,t)$ denote the degrees of membership and non-membership of $(s,t) \in S \times T$ in the relation $\tilde{R}$.

***Definition* 2.2.3:** Let $\tilde{P} = \langle \varsigma_{\tilde{P}}, \xi_{\tilde{P}} \rangle$ and $\tilde{Q} = \langle \varsigma_{\tilde{Q}}, \xi_{\tilde{Q}} \rangle$ be two intuitionistic fuzzy sets defined in $U$. Further suppose $\tilde{P} = \langle \mu_{\tilde{P}}, v_{\tilde{P}} \rangle$ is an intuitionistic fuzzy relation on $U$. Then $\tilde{P} = \langle \mu_{\tilde{P}}, v_{\tilde{P}} \rangle$ is said to be an intuitionistic fuzzy relation on $\tilde{Q} = \langle \varsigma_{\tilde{Q}}, \xi_{\tilde{Q}} \rangle$ if $\varsigma_{\tilde{P}}(u,v) \leq \min(\varsigma_{\tilde{Q}}(u), \varsigma_{\tilde{Q}}(v))$ and $\xi_{\tilde{P}}(u,v) \geq \max(\xi_{\tilde{Q}}(u), \xi_{\tilde{Q}}(v)), \forall u,v \in U.$

An IFR $\tilde{P}$ on $U$ is called symmetric if $\varsigma_{\tilde{P}}(u,v) = \varsigma_{\tilde{P}}(v,u)$ and $\xi_{\tilde{P}}(u,v) = \xi_{\tilde{P}}(v,u) \quad \forall u,v \in U.$

***Definition* 2.2.4 (*Intuitionistic Fuzzy Graph*)** [39]: An intuitionistic fuzzy graph (IFG) with underlying set $V$ is a pair $G^* = (\tilde{P}, \tilde{Q})$, where $\tilde{P} = \langle \varsigma_{\tilde{P}}, \xi_{\tilde{P}} \rangle$ is an IFS in $V$ with $0 \leq \varsigma_{\tilde{P}}(m) + \xi_{\tilde{P}}(m) \leq 1 \quad \forall m \in V$ and $\tilde{Q} = \langle \varsigma_{\tilde{Q}}, \xi_{\tilde{Q}} \rangle$ is an IFS in $E \subseteq V \times V$ such that

$$\varsigma_{\tilde{Q}}(u,v) \leq \min(\varsigma_{\tilde{Q}}(u), \varsigma_{\tilde{Q}}(v)) \text{ and } \xi_{\tilde{Q}}(u,v) \leq \max(\xi_{\tilde{P}}(u), \xi_{\tilde{P}}(v)), \tag{7}$$

and $\quad 0 \leq \varsigma_{\tilde{Q}}(u,v) + \xi_{\tilde{Q}}(u,v) \leq 1 \quad \forall (m,n) \in E. \tag{8}$

Here $\tilde{P}$ and $\tilde{Q}$ represent the intuitionistic fuzzy vertex set of $G^*$ and the intuitionistic fuzzy edge set of $G^*$, respectively. Note that here $\tilde{Q}$ is a symmetric intuitionistic fuzzy relation on $\tilde{P}$.

***Definition* 2.2.5 (*Strong Intuitionistic Fuzzy Graph*)** [21]: An IFG $G^* = (\tilde{P}, \tilde{Q})$ is said to be a strong intuitionistic fuzzy graph if

$$\varsigma_{\tilde{Q}}(u,v) = \min(\varsigma_{\tilde{P}}(u), \xi_{\tilde{P}}(v)) \text{ and } \xi_{\tilde{Q}}(u,v) = \max(\xi_{\tilde{P}}(u), \xi_{\tilde{P}}(v)) \quad \forall uv \in E. \tag{9}$$

***Definition* 2.2.6 (*Complete Intuitionistic Fuzzy Graph*)** [21]: An IFG $G^* = (\tilde{P}, \tilde{Q})$ is said to be a complete intuitionistic fuzzy graph if

$$\varsigma_{\tilde{Q}}(u,v) = \min(\varsigma_{\tilde{P}}(u), \xi_{\tilde{P}}(v)) \text{ and } \xi_{\tilde{Q}}(u,v) = \max(\xi_{\tilde{P}}(u), \xi_{\tilde{P}}(v)) \ \forall u,v \in V. \quad (10)$$

## 2.3 *Pythagorean Fuzzy Set: Basic results*

The notion of Pythagorean fuzzy sets (PFSs) was introduced by Yager [50, 52] in 2013 to model the uncertain information in highly complex realistic problems that intuitionistic fuzzy sets cannot capture.

***Definition* 2.3.1 (*Pythagorean Fuzzy Set*)** [50, 52]: A Pythagorean fuzzy set $P$ defined in a finite universe of discourse $U = (u_1, u_2, \ldots, u_n)$ is given by

$$P = \{\langle x, \varsigma_P(u), \xi_P(u)\rangle | u \in U\}, \quad (11)$$

where

$$\varsigma_P(u): U \to [0,1], \xi_P(u): U \to [0,1] \text{ and } 0 \leq \varsigma_P^2(u) + \xi_P^2(u) \leq 1 \ \forall u \in U, \quad (12)$$

where the numbers $\varsigma_P(u)$ and $\xi_P(u)$ represent the degree of membership and the degree of non-membership, respectively, of $u \in U$ in $P$.

Additionally, for each Pythagorean fuzzy set $P$ in $U$, the hesitation degree or the Pythagorean index can be defined as

$$\tau_P(u) = \sqrt{1 - \varsigma_P^2(u) - \xi_P^2(u)} \ \forall u \in U. \quad (13)$$

Throughout this paper, the set of all Pythagorean fuzzy sets in $U$ will be represented by $PFS(U)$.

***Definition* 2.3.2 (*Set Operations on Pythagorean Fuzzy Set*)** [55]: Let $P, Q \in PFS(U)$ given by

$$P = \{\langle u, \varsigma_P(u), \xi_P(u)\rangle | u \in U\} \text{ and } Q = \{\langle u, \varsigma_Q(u), \xi_Q(u)\rangle | u \in U\}, \quad (14)$$

then some set operations can be defined as follows:

    i.    $P \subseteq Q$ if and only if $\varsigma_P(u) \leq \varsigma_Q(u)$ and $\xi_P(u) \geq \xi_Q(u) \ \forall u \in U$;

    ii.   $P = Q$ if and only if $P \subseteq Q$ and $Q \subseteq P$;

    iii.  $P^C = \{\langle u, \xi_P(u), \varsigma_P(u)\rangle | u \in U\}$;

    iv.  $P \cup Q = \{\langle u, \max(\varsigma_P(u), \varsigma_Q(u)), \min(\xi_P(u), \xi_Q(u))\rangle | u \in U\}$;

**v.** $P \cap Q = \{\langle u, \min(\varsigma_P(u), \varsigma_Q(u)), \max(\xi_P(u), \xi_Q(u))\rangle \mid u \in U\};$

In the next section, we propose the concept of Pythagorean fuzzy relation (PFR) and Pythagorean fuzzy graphs (PFGs).

## 3. Pythagorean Fuzzy Graphs

The Pythagorean fuzzy graphs are a new generalization of Euler's graph theory that represents complex graphical problems more appropriately. This approach has some advantages over intuitionistic fuzzy graphs in which the sum of the degrees of membership and nonmembership for any vertex or edge should be less than or equal to 1.

We start with the following formal definition of a Pythagorean fuzzy relation.

***Definition* 3.1 (*Pythagorean Fuzzy Relation*):** Let $S$ and $T$ be two arbitrary finite non-empty sets. A Pythagorean fuzzy relation $R(S \rightarrow T)$ in the universe $S \times T$ is a Pythagorean fuzzy set of the form:

$$R = \{\langle (s,t), \varsigma_R(s,t), \xi_R(s,t)\rangle \mid s \in S, t \in T\}, \tag{15}$$

where
$$\varsigma_R(s,t): S \times T \rightarrow [0,1] \text{ and } \xi_R(s,t): S \times T \rightarrow [0,1], \tag{16}$$

and
$$0 \leq \varsigma_R^2(s,t) + \xi_R^2(s,t) \leq 1 \quad \forall (s,t) \in S \times T. \tag{17}$$

Here the numbers $\varsigma_R(s,t)$ and $\xi_R(s,t)$ denote the degrees of membership and non-membership of the Pythagorean fuzzy relation $R$.

***Definition* 3.2:** Let $P = \langle \varsigma_P, \xi_P \rangle$ and $Q = \langle \varsigma_Q, \xi_Q \rangle$ be two Pythagorean fuzzy sets defined in $U$. Further suppose $P = \langle \varsigma_P, \xi_P \rangle$ is a Pythagorean fuzzy relation on $U$. Then $P = \langle \varsigma_P, \xi_P \rangle$ is called a PFR on $Q = \langle \varsigma_Q, \xi_Q \rangle$ if

$$\varsigma_P(u,v) \leq \min(\varsigma_Q(u), \varsigma_Q(v)) \text{ and } \xi_A(u,v) \geq \max(\xi_B(u), \xi_B(v)) \ \forall u, v \in U. \tag{18}$$

***Definition* 3.3 (*Symmetric Pythagorean Fuzzy Relation*):** A Pythagorean fuzzy relation $R$ on $U$ is called *symmetric* if

$$\varsigma_R(u,v) = \varsigma_R(v,u) \text{ and } \xi_R(u,v) = \xi_R(v,u) \ \forall u, v \in U. \tag{19}$$

***Definition* 3.4 (*Pythagorean Fuzzy Graph*):** A Pythagorean fuzzy graph with the underlying set $V$ is a pair $G^{**} = (P, Q)$, where $P = \langle \varsigma_P, \xi_P \rangle$ is a Pythagorean fuzzy set in $V$ with $0 \leq \varsigma_P^2(u) + \xi_P^2(u) \leq 1 \ \forall u \in V$ and $Q = \langle \varsigma_Q, \xi_Q \rangle$ is a Pythagorean fuzzy set in $E \subseteq V \times V$ such that

$$\varsigma_Q(u,v) \leq \min(\varsigma_P(u), \varsigma_P(v)) \text{ and } \xi_B(u,v) \leq \max(\xi_P(u), \xi_P(v)), \tag{20}$$

and
$$0 \leq \varsigma_Q^2(u,v) + \xi_Q^2(u,v) \leq 1 \ \forall (u,v) \in E. \tag{21}$$

Then, $P$ and $Q$ are the Pythagorean fuzzy vertex set and the Pythagorean fuzzy edge set of $G^{**}$, respectively. Note that here $Q$ is a symmetric Pythagorean fuzzy relation on $P$. For convenience, we use the notation $uv$ for an element of $E$.

Throughout this paper, $G$, $G^*$ and $G^{**}$ represent, respectively, a crisp graph, an intuitionistic fuzzy graph, and a Pythagorean fuzzy graph.

***Remark* 3.1**: **(i)** When $\varsigma_Q(uv) = 0$ and $\xi_Q(uv) = 0$ for some $u$ and $v$, then there is no edge between $u$ and $v$.

**(ii)** When either one of the following conditions is true, then there is an edge between $u$ and $v$:

**(a)** $\varsigma_Q(uv) = 0$ and $\xi_Q(uv) > 0$,

**(b)** $\varsigma_Q(uv) > 0$ and $\xi_Q(uv) = 0$,

**(c)** $\varsigma_Q(uv) > 0$ and $\xi_Q(uv) > 0$.

***Example* 3.1:** Let $G = (V, E)$ be a graph such that $V = \{a, b, c, d\}$ and $E = \{ab, bc, cd, da\} \subseteq V \times V$. Let $P = \langle \varsigma_P, \xi_P \rangle$ be a Pythagorean fuzzy set in $V$ and $Q = \langle \varsigma_Q, \xi_Q \rangle$ be a Pythagorean fuzzy set in $E$ defined by

$$\left. \begin{array}{l} P = \langle (a, 0.5, 0.7), (b, 0.8, 0.3), (c, 0.6, 0.5), (d, 0.4, 0.4) \rangle, \\ Q = \langle (ab, 0.4, 0.7), (bc, 0.5, 0.45), (cd, 0.3, 0.5), (da, 0.4, 0.6) \rangle \end{array} \right\}. \tag{22}$$

Then the ordered pair $G^{**} = (P, Q)$ is a Pythagorean fuzzy graph of $G$.

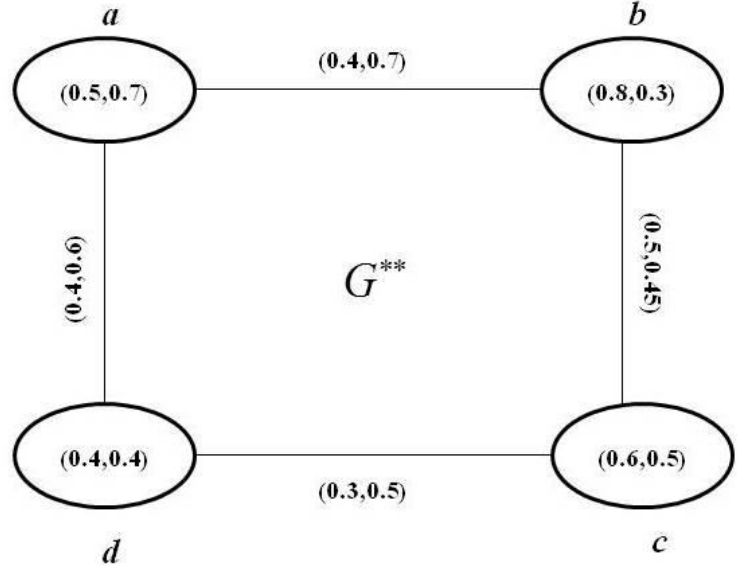

We now turn to study some basic operations on Pythagorean fuzzy graphs:

***Definition 3.5 (Cartesian Product of Pythagorean Fuzzy Graphs):*** Let $G_1^{**} = (P_1, Q_1)$ and $G_2^{**} = (P_2, Q_2)$ be two Pythagorean fuzzy graphs of the graphs $G_1 = (V_1, E_1)$ and $G_2 = (V_2, E_2)$, respectively, where $P_1 = \langle \varsigma_{P_1}, \xi_{P_1} \rangle$ and $P_2 = \langle \varsigma_{P_2}, \xi_{P_2} \rangle$ are Pythagorean fuzzy sets, correspondingly, defined in $V_1$ and $V_2$; $Q_1 = \langle \varsigma_{Q_1}, \xi_{Q_1} \rangle$ and $Q_2 = \langle \varsigma_{Q_2}, \xi_{Q_2} \rangle$ are Pythagorean fuzzy sets in $E_1 \subseteq V_1 \times V_1$ and $E_2 \subseteq V_2 \times V_2$, respectively.

The Cartesian product of graphs $G_1^{**}$ and $G_2^{**}$, denoted by $G_1^{**} \times G_2^{**} = (P_1 \times P_2, Q_1 \times Q_2)$, is defined as follows

i. $\varsigma_{P_1 \times P_2}(u_1, u_2) = \min(\varsigma_{P_1}(u_1), \varsigma_{P_2}(u_2))$ and

$\xi_{P_1 \times P_2}(u_1, u_2) = \max(\xi_{P_1}(u_1), \xi_{P_2}(u_2)) \quad \forall (u_1, u_2) \in V_1 \times V_2 = V.$

ii. $\varsigma_{Q_1 \times Q_2}((u, u_2)(u, v_2)) = \min(\varsigma_{P_1}(u), \varsigma_{Q_2}(u_2 v_2))$ and

$\xi_{Q_1 \times Q_2}((u, u_2)(u, v_2)) = \max(\xi_{P_1}(u), \xi_{Q_2}(u_2 v_2)) \quad \forall u \in V_1, u_2 v_2 \in E_2.$

iii. $\varsigma_{Q_1 \times Q_2}((u_1, w)(v_1, w)) = \min(\varsigma_{Q_1}(u_1 v_1), \varsigma_{P_2}(w))$ and

$\xi_{Q_1 \times Q_2}((u_1, w)(v_1, w)) = \max(\xi_{Q_1}(u_1 v_1), \xi_{P_2}(w)), \quad \forall w \in V_2, u_1 v_1 \in E_1.$

***Proposition* 3.1:** The Cartesian product of two Pythagorean fuzzy graphs is a Pythagorean fuzzy graph.

***Proof*:** Let $t \in V_1$, $u_2 v_2 \in E_2$, then we have

$$\varsigma_{Q_1 \times Q_2}((t,u_2)(t,v_2)) = \min(\varsigma_{P_1}(t), \varsigma_{Q_2}(u_2 v_2))$$
$$\leq \min(\varsigma_{P_1}(t), \min(\varsigma_{P_2}(u_2), \varsigma_{P_2}(v_2)))$$
$$= \min(\min(\varsigma_{P_1}(t), \varsigma_{P_2}(u_2)), \min(\varsigma_{P_1}(t), \varsigma_{P_2}(v_2)))$$
$$= \min(\varsigma_{P_1 \times P_2}(t,u_2), \varsigma_{P_1 \times P_2}(t,v_2)), \tag{23}$$

$$\xi_{Q_1 \times Q_2}((t,u_2)(t,v_2)) = \max(\xi_{P_1}(t), \xi_{Q_2}(u_2 v_2))$$
$$\leq \max(\xi_{P_1}(t), \max(\xi_{P_2}(u_2), \xi_{P_2}(v_2)))$$
$$= \max(\max(\xi_{A_1}(t), \xi_{A_2}(u_2)), \max(\xi_{A_1}(t), \xi_{A_2}(v_2)))$$
$$= \max(\xi_{P_1 \times P_2}(t,u_2), \xi_{P_1 \times P_2}(t,v_2)). \tag{24}$$

Again, let $w \in V_2$, $u_1 v_1 \in E_1$, then we have

$$\varsigma_{Q_1 \times Q_2}((u_1,w)(v_1,w)) = \min(\varsigma_{Q_1}(u_1 v_1), \varsigma_{P_2}(w))$$
$$\leq \min(\min(\varsigma_{P_1}(u_1), \varsigma_{P_1}(v_1)), \varsigma_{P_2}(w))$$
$$= \min(\min(\varsigma_{P_1}(u_1), \varsigma_{P_2}(w)), \min(\varsigma_{P_1}(v_1), \varsigma_{P_2}(w)))$$
$$= \min(\varsigma_{P_1 \times P_2}(u_1,w), \varsigma_{P_1 \times P_2}(v_1,w)), \tag{25}$$

$$\xi_{Q_1 \times Q_2}((u_1,w)(v_1,w)) = \max(\xi_{Q_1}(u_1 v_1), \xi_{P_2}(w))$$
$$\leq \max(\max(\xi_{P_1}(u_1), \xi_{P_1}(v_1)), \xi_{P_2}(w))$$
$$= \max(\max(\xi_{P_1}(u_1), \xi_{P_2}(w)), \max(\xi_{P_1}(v_1), \xi_{P_2}(w)))$$
$$= \max(\xi_{P_1 \times P_2}(u_1,w), \xi_{P_1 \times P_2}(v_1,w)). \tag{26}$$

This completes the proof. □

***Example* 3.2:** Consider $G_1 = (V_1, E_1)$ and $G_2 = (V_2, E_2)$ are two graphs such that $V_1 = \{a,b\}$, $V_2 = \{c,d\}$, $E_1 = \{ab\}$ and $E_2 = \{cd\}$. Let $G_1^{**} = (P_1, Q_1)$ and $G_2^{**} = (P_2, Q_2)$ be Pythagorean fuzzy graphs in $G_1$ and $G_2$, respectively, where

$$P_1 = \langle (a, 0.6, 0.3), (b, 0.5, 0.7) \rangle, Q_1 = \langle (ab, 0.5, 0.7) \rangle. \tag{27}$$

$$P_2 = \langle (c, 0.7, 0.5), (d, 0.5, 0.8) \rangle, Q_2 = \langle (cd, 0.4, 0.65) \rangle. \tag{28}$$

It is easy to verify that $G_1^{**} \times G_2^{**} = (P_1 \times P_2, Q_1 \times Q_2)$ is a PFG of $G_1 \times G_2$.

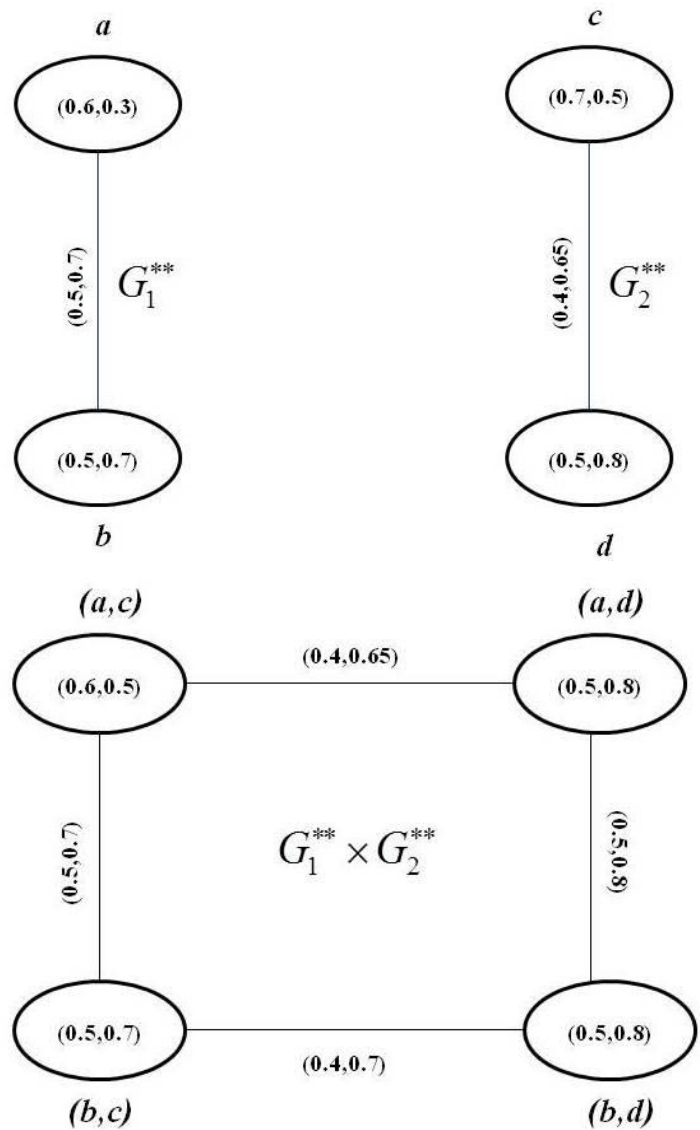

***Definition* 3.6 (*Composition of Pythagorean Fuzzy Graphs*):** Let $G_1^{**} = (P_1, Q_1)$ and $G_2^{**} = (P_2, Q_2)$ be two Pythagorean fuzzy graphs of the graphs $G_1 = (V_1, E_1)$ and $G_2 = (V_2, E_2)$, respectively, where $P_1 = \langle \varsigma_{P_1}, \xi_{P_1} \rangle$ and $P_2 = \langle \varsigma_{P_2}, \xi_{P_2} \rangle$ are Pythagorean fuzzy sets, correspondingly, defined in $V_1$ and $V_2$; $Q_1 = \langle \varsigma_{Q_1}, \xi_{Q_1} \rangle$ and $Q_2 = \langle \varsigma_{Q_2}, \xi_{Q_2} \rangle$ are Pythagorean fuzzy sets in $E_1 \subseteq V_1 \times V_1$ and $E_2 \subseteq V_2 \times V_2$, respectively.

The composition of graphs $G_1^{**}$ and $G_2^{**}$, denoted by $G_1^{**}[G_2^{**}] = (P_1 \circ P_2, Q_1 \circ Q_2)$, is defined as follows

i. $\varsigma_{P_1 \circ P_2}(u_1, u_2) = \min(\varsigma_{P_1}(u_1), \varsigma_{P_2}(u_2))$ and

$\xi_{P_1 \circ P_2}(u_1, u_2) = \max(\xi_{P_1}(u_1), \xi_{P_2}(u_2)) \quad \forall (u_1, u_2) \in V_1 \times V_2 = V$.

ii. $\varsigma_{Q_1 \circ Q_2}((t, u_2)(t, v_2)) = \min(\varsigma_{P_1}(t), \varsigma_{Q_2}(u_2 v_2))$ and

$\xi_{Q_1 \circ Q_2}((t, u_2)(t, u_2)) = \max(\xi_{P_1}(t), \xi_{Q_2}(u_2 v_2)) \quad \forall\, t \in V_1, u_2 v_2 \in E_2$.

iii. $\varsigma_{Q_1 \circ Q_2}((u_1, w)(v_1, w)) = \min(\varsigma_{Q_1}(u_1 v_1), \varsigma_{P_2}(w))$ and

$\xi_{Q_1 \circ Q_2}((u_1, w)(v_1, w)) = \max(\xi_{Q_1}(u_1 v_1), \xi_{P_2}(w)) \quad \forall\, w \in V_2, x_1 y_1 \in E_1$.

iv. $\varsigma_{Q_1 \circ Q_2}((u_1, u_2)(v_1, v_2)) = \min(\varsigma_{P_2}(u_2), \varsigma_{P_2}(v_2), \varsigma_{Q_1}(u_1 v_1))$ and

$\xi_{Q_1 \circ Q_2}((u_1, u_2)(v_1, v_2)) = \max(\xi_{P_2}(u_2), \xi_{P_2}(v_2), \xi_{Q_1}(u_1 v_1)) \quad \forall\, (u_1, u_2)(v_1, v_2) \in E^0 - E$.

where $E^0 = E \cup \{(u_1, u_2)(v_1, v_2) \mid u_1 v_1 \in E_1, u_2 \neq v_2\}$ and $E$ is defined as in the case for $G_1^{**} \times G_2^{**}$.

***Proposition* 3.2:** The composition of two Pythagorean fuzzy graphs is a Pythagorean fuzzy graph.

***Proof*:** Let $t \in V_1$, $u_2 v_2 \in E_2$, then we have

$$\varsigma_{Q_1 \circ Q_2}((t, u_2)(t, v_2)) = \min(\varsigma_{P_1}(t), \varsigma_{Q_2}(u_2 v_2))$$

$$\leq \min(\varsigma_{P_1}(t), \min(\varsigma_{P_2}(u_2), \varsigma_{P_2}(v_2)))$$

$$= \min(\min(\varsigma_{P_1}(t), \varsigma_{P_2}(u_2)), \min(\varsigma_{P_1}(t), \varsigma_{P_2}(v_2)))$$

$$= \min(\varsigma_{P_1 \circ P_2}(t, u_2), \varsigma_{P_1 \circ P_2}(t, v_2)), \tag{29}$$

$$\xi_{Q_1 \circ Q_2}\big((t,u_2)(t,v_2)\big) = \max\big(\xi_{P_1}(t), \xi_{Q_2}(u_2 v_2)\big)$$

$$\leq \max\big(\xi_{P_1}(t), \max\big(\xi_{P_2}(u_2), \xi_{P_2}(v_2)\big)\big)$$

$$= \max\big(\max\big(\xi_{P_1}(t), \xi_{P_2}(u_2)\big), \max\big(\xi_{P_1}(t), \xi_{P_2}(v_2)\big)\big)$$

$$= \max\big(\xi_{P_1 \circ P_2}(t, u_2), \xi_{P_1 \circ P_2}(t, v_2)\big). \tag{30}$$

Again, let $w \in V_2$, $u_1 v_1 \in E_1$, then we have

$$\varsigma_{Q_1 \circ Q_2}\big((u_1, w)(v_1, w)\big) = \min\big(\varsigma_{Q_1}(u_1 v_1), \varsigma_{P_2}(w)\big)$$

$$\leq \min\big(\min\big(\varsigma_{P_1}(u_1), \varsigma_{P_1}(v_1)\big), \mu_{P_2}(w)\big)$$

$$= \min\big(\min\big(\varsigma_{P_1}(u_1), \varsigma_{P_2}(w)\big), \min\big(\varsigma_{P_1}(v_1), \varsigma_{P_2}(w)\big)\big)$$

$$= \min\big(\varsigma_{P_1 \circ P_2}(u_1, w), \varsigma_{P_1 \circ P_2}(v_1, w)\big), \tag{31}$$

$$\xi_{Q_1 \circ Q_2}\big((u_1, w)(v_1, w)\big) = \max\big(\xi_{Q_1}(u_1 v_1), \xi_{P_2}(w)\big)$$

$$\leq \max\big(\max\big(\xi_{P_1}(u_1), \xi_{P_1}(v_1)\big), \xi_{P_2}(w)\big)$$

$$= \max\big(\max\big(\xi_{P_1}(u_1), \xi_{P_2}(w)\big), \max\big(\xi_{P_1}(v_1), \xi_{P_2}(w)\big)\big)$$

$$= \max\big(\xi_{P_1 \circ P_2}(u_1, w), \xi_{P_1 \circ P_2}(v_1, w)\big). \tag{32}$$

Further, let $(u_1, u_2)(v_1, v_2) \in E^0 - E$, so $u_1 v_1 \in E_1$, $u_2 \neq v_2$. Then it follows that

$$\varsigma_{Q_1 \circ Q_2}\big((u_1, u_2)(v_1, v_2)\big) = \min\big(\varsigma_{P_2}(u_2), \varsigma_{P_2}(v_2), \varsigma_{Q_1}(u_1 v_1)\big)$$

$$\leq \min\big(\varsigma_{P_2}(u_2), \varsigma_{P_2}(v_2), \min\big(\varsigma_{P_1}(u_1), \varsigma_{P_1}(v_1)\big)\big)$$

$$= \min\big(\min\big(\varsigma_{P_1}(u_1), \varsigma_{P_2}(u_2)\big), \min\big(\varsigma_{P_1}(v_1), \varsigma_{P_2}(v_2)\big)\big)$$

$$= \min\big(\varsigma_{P_1 \circ P_2}(u_1, u_2), \varsigma_{P_1 \circ P_2}(v_1, v_2)\big), \tag{33}$$

$$\xi_{Q_1 \circ Q_2}\big((u_1, u_2)(v_1, v_2)\big) = \max\big(\xi_{P_2}(u_2), \xi_{P_2}(v_2), \xi_{Q_1}(u_1 v_1)\big)$$

$$\leq \max\big(\xi_{P_2}(u_2), \xi_{P_2}(v_2), \max\big(\xi_{P_1}(u_1), \xi_{P_1}(v_1)\big)\big)$$

$$= \max\big(\max\big(\xi_{P_1}(u_1), \xi_{P_2}(u_2)\big), \min\big(\xi_{P_1}(v_1), \xi_{P_2}(v_2)\big)\big)$$

$$= \max\big(\xi_{P_1 \circ P_2}(u_1, u_2), \xi_{P_1 \circ P_2}(v_1, v_2)\big). \tag{34}$$

This completes the proof. □

**Example 3.3:** Let $G_1 = (V_1, E_1)$ and $G_2 = (V_2, E_2)$ be graphs such that $V_1 = \{a, b\}, V_2 = \{c, d\}$, $E_1 = \{ab\}$ and $E_2 = \{cd\}$. Consider two Pythagorean fuzzy graphs $G_1^{**} = (P_1, Q_1)$ and $G_2^{**} = (P_2, Q_2)$, where

$$P_1 = \langle (a, 0.6, 0.3), (b, 0.5, 0.7) \rangle, Q_1 = \langle (ab, 0.5, 0.7) \rangle . \tag{35}$$

$$P_2 = \langle (c, 0.7, 0.5), (d, 0.5, 0.8) \rangle, Q_2 = \langle (cd, 0.4, 0.65) \rangle . \tag{36}$$

By simple computation, it is easy to verify that the graph $G_1^{**}[G_2^{**}]$ is a PFG of $G_1[G_2]$.

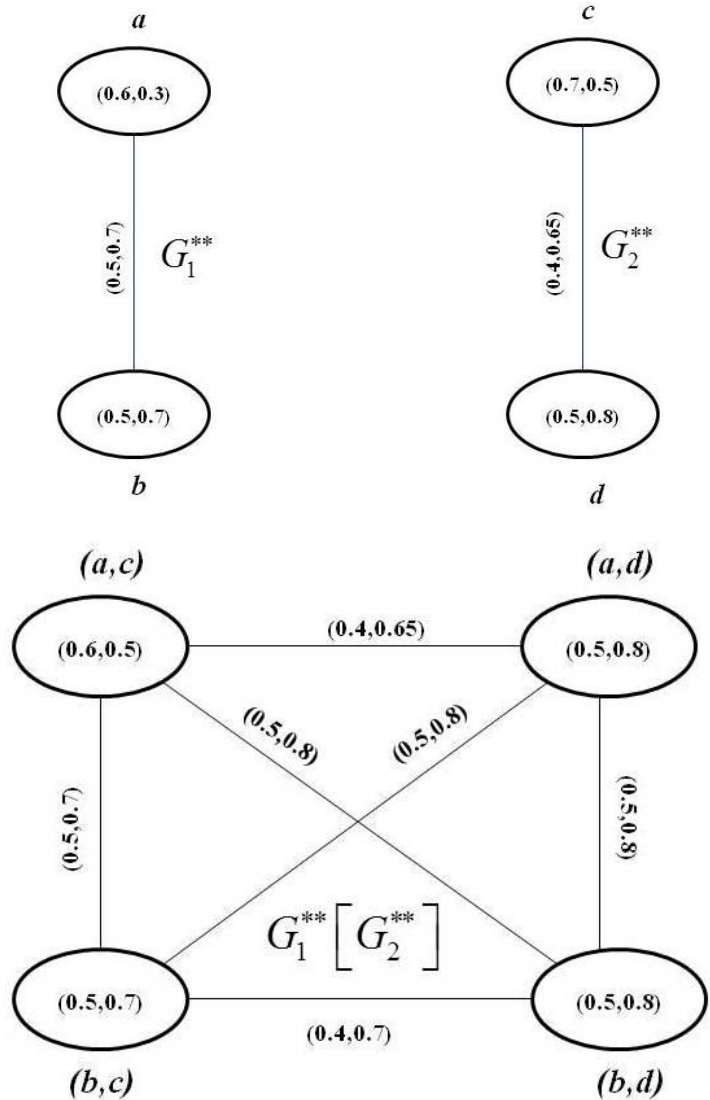

***Definition* 3.7 (*Union of Pythagorean Fuzzy Graphs*):** Let $G_1^{**} = (P_1, Q_1)$ and $G_2^{**} = (P_2, Q_2)$ be two Pythagorean fuzzy graphs of the graphs $G_1 = (V_1, E_1)$ and $G_2 = (V_2, E_2)$ respectively, where $P_1 = \langle \varsigma_{P_1}, \xi_{P_1} \rangle$ and $P_2 = \langle \varsigma_{P_2}, \xi_{P_2} \rangle$ are Pythagorean fuzzy sets, correspondingly, defined in $V_1$ and $V_2$; $Q_1 = \langle \varsigma_{Q_1}, \xi_{Q_1} \rangle$ and $Q_2 = \langle \varsigma_{Q_2}, \xi_{Q_2} \rangle$ are Pythagorean fuzzy sets in $E_1 \subseteq V_1 \times V_1$ and $E_2 \subseteq V_2 \times V_2$, respectively.

The union of two Pythagorean fuzzy graphs $G_1^{**}$ and $G_2^{**}$, denoted by $G_1^{**} \cup G_2^{**} = (P_1 \cup P_2, Q_1 \cup Q_2)$, is defined as follows

i. $\begin{cases} \varsigma_{P_1 \cup P_2}(u) = \varsigma_{P_1}(u) \text{ if } u \in V_1 \cap \overline{V_2} \\ \varsigma_{P_1 \cup P_2}(u) = \varsigma_{P_2}(u) \text{ if } u \in V_2 \cap \overline{V_1} \\ \varsigma_{P_1 \cup P_2}(u) = \max\left(\varsigma_{P_1}(u), \varsigma_{P_2}(u)\right) \text{ if } u \in V_1 \cap V_2 \end{cases}$.

ii. $\begin{cases} \xi_{P_1 \cup P_2}(u) = \xi_{A_1}(u) \text{ if } u \in V_1 \cap \overline{V_2} \\ \xi_{P_1 \cup P_2}(u) = \xi_{A_2}(u) \text{ if } u \in V_2 \cap \overline{V_1} \\ \xi_{P_1 \cup P_2}(u) = \min\left(\xi_{A_1}(u), \xi_{A_2}(u)\right) \text{ if } u \in V_1 \cap V_2 \end{cases}$.

iii. $\begin{cases} \varsigma_{Q_1 \cup Q_2}(uv) = \varsigma_{Q_1}(uv) \text{ if } uv \in E_1 \cap \overline{E_2} \\ \varsigma_{Q_1 \cup Q_2}(uv) = \varsigma_{Q_2}(uv) \text{ if } uv \in E_2 \cap \overline{E_1} \\ \varsigma_{Q_1 \cup Q_2}(uv) = \max\left(\varsigma_{Q_1}(uv), \varsigma_{Q_2}(uv)\right) \text{ if } uv \in E_1 \cap E_2 \end{cases}$

iv. $\begin{cases} \xi_{Q_1 \cup Q_2}(uv) = \xi_{Q_1}(uv) \text{ if } uv \in E_1 \cap \overline{E_2} \\ \xi_{Q_1 \cup Q_2}(uv) = \xi_{Q_2}(uv) \text{ if } uv \in E_2 \cap \overline{E_1} \\ \xi_{Q_1 \cup Q_2}(uv) = \min\left(\xi_{Q_1}(uv), \xi_{Q_2}(uv)\right) \text{ if } uv \in E_1 \cap E_2. \end{cases}$

***Proposition* 3.3:** The union of two Pythagorean fuzzy graphs is a Pythagorean fuzzy graph.

***Proof*:** Let $uv \in E_1 \cap E_2$, then we have

$\varsigma_{Q_1 \cup Q_2}(uv) = \max\left(\varsigma_{Q_1}(uv), \varsigma_{Q_2}(uv)\right)$

$\leq \max\left(\min\left(\varsigma_{P_1}(u), \varsigma_{P_1}(v)\right), \min\left(\varsigma_{P_2}(u), \varsigma_{P_2}(v)\right)\right)$

$\leq \min\left(\max\left(\varsigma_{P_1}(u), \varsigma_{P_2}(u)\right), \max\left(\varsigma_{P_1}(v), \varsigma_{P_2}(v)\right)\right)$

$$= \min\left(\varsigma_{P_1 \cup P_2}(u), \varsigma_{P_1 \cup P_2}(v)\right), \tag{37}$$

$$\xi_{Q_1 \cup Q_2}(uv) = \min\left(\xi_{Q_1}(uv), \xi_{Q_2}(uv)\right)$$

$$\leq \min\left(\max\left(\nu_{P_1}(u), \nu_{P_1}(v)\right), \max\left(\nu_{P_2}(u), \nu_{P_2}(v)\right)\right)$$

$$= \max\left(\min\left(\xi_{P_1}(u), \xi_{P_2}(u)\right), \min\left(\xi_{P_1}(v), \xi_{P_2}(v)\right)\right)$$

$$= \max\left(\xi_{P_1 \cup P_2}(u), \xi_{P_1 \cup P_2}(v)\right). \tag{38}$$

Again, let $uv \in E_1 \cap \bar{E}_2$, then

$$\varsigma_{Q_1 \cup Q_2}(uv) = \varsigma_{Q_1}(uv)$$

$$\leq \min\left(\varsigma_{P_1}(u), \varsigma_{P_1}(v)\right)$$

$$= \min\left(\mu_{P_1 \cup P_2}(u), \mu_{P_1 \cup P_2}(v)\right), \quad \text{if } u, v \in V_1 \cap \bar{V}_2$$

$$\leq \min\left(\varsigma_{P_1 \cup P_2}(u), \max\left(\varsigma_{P_1}(v), \varsigma_{P_2}(v)\right)\right)$$

$$= \min\left(\varsigma_{P_1 \cup P_2}(u), \varsigma_{P_1 \cup P_2}(v)\right), \quad \text{if } u \in V_1 \cap \bar{V}_2, v \in V_1 \cap V_2$$

$$\leq \min\left(\max\left(\varsigma_{P_1}(u), \varsigma_{P_2}(u)\right), \max\left(\varsigma_{P_1}(v), \varsigma_{P_2}(v)\right)\right)$$

$$= \min\left(\varsigma_{P_1 \cup P_2}(u), \varsigma_{P_1 \cup P_2}(v)\right), \quad \text{if } u, v \in V_1 \cap V_2. \tag{39}$$

$$\xi_{Q_1 \cup Q_2}(uv) = \xi_{Q_1}(uv)$$

$$\leq \max\left(\xi_{P_1}(u), \xi_{P_1}(v)\right)$$

$$= \max\left(\xi_{P_1 \cup P_2}(u), \xi_{P_1 \cup P_2}(v)\right), \quad \text{if } u, v \in V_1 \cap \bar{V}_2$$

$$\leq \max\left(\xi_{P_1 \cup P_2}(u), \min\left(\xi_{P_1}(v), \xi_{P_2}(v)\right)\right)$$

$$= \max\left(\xi_{P_1 \cup P_2}(u), \xi_{P_1 \cup P_2}(v)\right), \quad \text{if } u \in V_1 \cap \bar{V}_2, v \in V_1 \cap V_2$$

$$\leq \max\left(\min\left(\xi_{P_1}(u), \xi_{P_2}(u)\right), \min\left(\xi_{P_1}(v), \xi_{P_2}(v)\right)\right)$$

$$= \max\left(\xi_{P_1 \cup P_2}(u), \xi_{P_1 \cup P_2}(v)\right), \text{ if } x, y \in V_1 \cap V_2. \tag{40}$$

Similarly if $uv \in \bar{E}_1 \cap E_2$, then

$$\varsigma_{Q_1 \cup Q_2}(uv) \leq \min\left(\varsigma_{P_1 \cup P_2}(u), \varsigma_{P_1 \cup P_2}(v)\right) \text{ and } \xi_{Q_1 \cup Q_2}(uv) \leq \max\left(\xi_{P_1 \cup P_2}(u), \xi_{P_1 \cup P_2}(v)\right). \tag{41}$$

This completes the proof. □

**Example 3.4:** Let $G_1 = (V_1, E_1)$ and $G_2 = (V_2, E_2)$ be graphs such that $V_1 = \{a, b, c, d, e\}$, $V_2 = \{a, b, c, d, f\}$, $E_1 = \{ab, bc, be, ce, ad, de\}$ and $E_2 = \{ab, bc, bf, bd, cf\}$. Consider two Pythagorean fuzzy graphs $G_1^{**} = (P_1, Q_1)$ and $G_2^{**} = (P_2, Q_2)$, where

$$\left. \begin{array}{l} P_1 = \langle (a, 0.3, 0.8), (b, 0.5, 0.6), (c, 0.3, 0.4), (d, 0.7, 0.2), (e, 0.6, 0.6) \rangle, \\ Q_1 = \langle (ab, 0.3, 0.7), (bc, 0.3, 0.6)(be, 0.5, 0.6), (ce, 0.2, 0.5), (de, 0.5, 0.6), (ad, 0.2, 0.8) \rangle \end{array} \right\}, \quad (42)$$

$$\left. \begin{array}{l} P_2 = \langle (a, 0.7, 0.1), (b, 0.4, 0.6), (c, 0.8, 0.2), (d, 0.2, 0.4), (f, 0.6, 0.7) \rangle, \\ Q_2 = \langle (ab, 0.4, 0.6), (bc, 0.3, 0.6)(bf, 0.4, 0.7), (bd, 0.2, 0.6), (cf, 0.5, 0.7) \rangle \end{array} \right\}. \quad (43)$$

Clearly, $G_1^{**} \cup G_2^{**}$ is a Pythagorean fuzzy graph of $G_1 \cup G_2$.

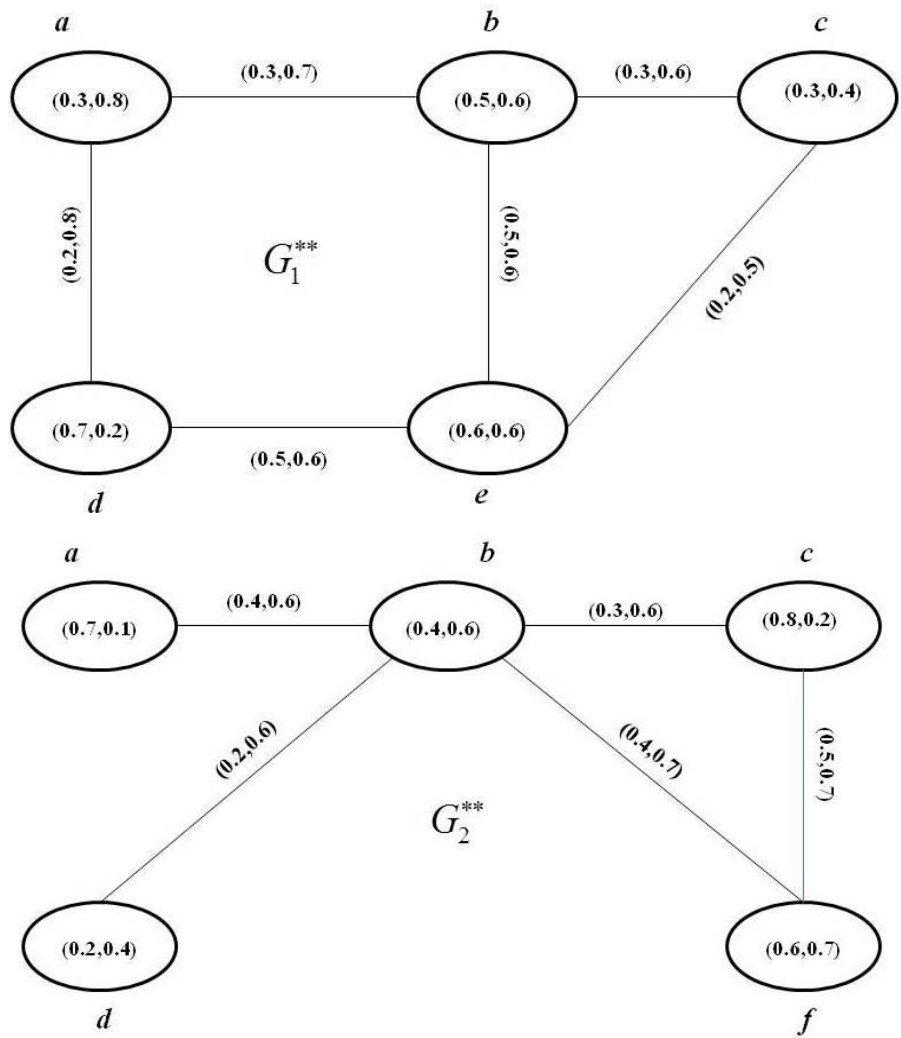

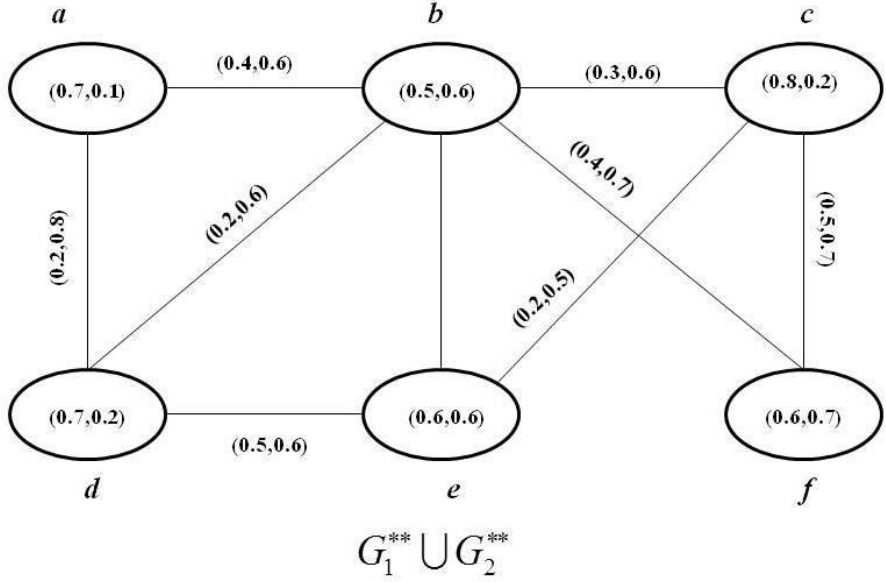

$$G_1^{**} \cup G_2^{**}$$

**Definition 3.8 (Join of Pythagorean Fuzzy Graphs):** Let $G_1^{**} = (P_1, Q_1)$ and $G_2^{**} = (P_2, Q_2)$ be two Pythagorean fuzzy graphs of the graphs $G_1 = (V_1, E_1)$ and $G_2 = (V_2, E_2)$ respectively, where $P_1 = \langle \varsigma_{P_1}, \xi_{P_1} \rangle$ and $P_2 = \langle \varsigma_{P_2}, \xi_{P_2} \rangle$ are Pythagorean fuzzy sets, correspondingly, defined in $V_1$ and $V_2$; $Q_1 = \langle \varsigma_{Q_1}, \xi_{Q_1} \rangle$ and $Q_2 = \langle \varsigma_{Q_2}, \xi_{Q_2} \rangle$ are Pythagorean fuzzy sets in $E_1 \subseteq V_1 \times V_1$ and $E_2 \subseteq V_2 \times V_2$, respectively.

The join of two Pythagorean fuzzy graphs $G_1^{**}$ and $G_2^{**}$, denoted by $G_1^{**} + G_2^{**} = (P_1 + P_2, Q_1 + Q_2)$, is defined as follows

i. $\begin{cases} \varsigma_{P_1+P_2}(u) = \varsigma_{P_1}(u) & \text{if } u \in V_1 \cap \overline{V_2}, \\ \varsigma_{P_1+P_2}(u) = \varsigma_{P_2}(u) & \text{if } u \in \overline{V_1} \cap V_2, \\ \varsigma_{P_1+P_2}(u) = \varsigma_{P_1 \cup P_2}(u) & \text{if } u \in V_1 \cap V_2. \end{cases}$

ii. $\begin{cases} \xi_{P_1+P_2}(u) = v_{P_1}(u) & \text{if } u \in V_1 \cap \overline{V_2}, \\ \xi_{P_1+P_2}(u) = v_{P_2}(u) & \text{if } u \in \overline{V_1} \cap V_2, \\ \xi_{P_1+P_2}(u) = v_{P_1 \cup P_2}(u) & \text{if } u \in V_1 \cap V_2. \end{cases}$

iii. $\begin{cases} \varsigma_{Q_1+Q_2}(uv) = \varsigma_{Q_1}(uv) & \text{if } uv \in E_1 \cap \overline{E_2}, \\ \varsigma_{Q_1+Q_2}(uv) = \varsigma_{Q_2}(uv) & \text{if } uv \in \overline{E_1} \cap E_2, \\ \varsigma_{Q_1+Q_2}(uv) = \varsigma_{Q_1 \cup Q_2}(uv) & \text{if } uv \in E_1 \cap E_2. \end{cases}$

**iv.**
$$\begin{cases} \xi_{Q_1+Q_2}(uv) = \xi_{Q_1}(uv) & \text{if } uv \in E_1 \cap \bar{E}_2, \\ \xi_{Q_1+Q_2}(uv) = \xi_{Q_2}(uv) & \text{if } uv \in \bar{E}_1 \cap E_2, \\ \xi_{Q_1+Q_2}(uv) = \xi_{Q_1 \cup Q_2}(uv) & \text{if } uv \in E_1 \cap E_2. \end{cases}$$

**v.**
$$\begin{cases} \varsigma_{Q_1+Q_2}(uv) = \min\left(\varsigma_{P_1}(u), \varsigma_{P_2}(v)\right) \\ \xi_{Q_1+Q_2}(uv) = \max\left(\xi_{P_1}(u), \xi_{P_2}(v)\right) \end{cases} \text{ if } uv \in E',$$

where $E'$ represents the set of all edges joining the nodes of $V_1$ and $V_2$.

***Proposition* 3.4:** The join of two Pythagorean fuzzy graphs is a Pythagorean fuzzy graph.

***Proof*:** Let $uv \in E'$, then we have

$$\varsigma_{Q_1+Q_2}(uv) = \min\left(\varsigma_{P_1}(u), \varsigma_{P_2}(v)\right)$$
$$\leq \min\left(\max\left(\varsigma_{P_1}(u), \varsigma_{P_2}(u)\right), \max\left(\varsigma_{P_1}(v), \varsigma_{P_2}(v)\right)\right)$$
$$= \min\left(\varsigma_{P_1+P_2}(u), \varsigma_{P_1+P_2}(v)\right), \tag{44}$$

$$\xi_{Q_1+Q_2}(uv) = \max\left(\xi_{P_1}(u), \xi_{P_2}(v)\right)$$
$$\leq \max\left(\min\left(\xi_{P_1}(u), \xi_{P_2}(u)\right), \min\left(\xi_{P_1}(v), \xi_{P_2}(v)\right)\right)$$
$$= \max\left(\xi_{P_1+P_2}(u), \xi_{P_1+P_2}(v)\right). \tag{45}$$

Again, let $uv \in E_1 \cup E_2$, the desired result directly obtains from Proposition 3.3.

This proves the proposition. □

***Example* 3.5:** Let $G_1 = (V_1, E_1)$ and $G_2 = (V_2, E_2)$ be graphs such that $V_1 = \{a, b\}, V_2 = \{c, d, e\}$, $E_1 = \{ab\}$ and $E_2 = \{cd, de\}$. Consider two Pythagorean fuzzy graphs $G_1^{**} = (P_1, Q_1)$ and $G_2^{**} = (P_2, Q_2)$, where

$$P_1 = \langle(a, 0.6, 0.3), (b, 0.5, 0.7)\rangle, Q_1 = \langle(ab, 0.5, 0.7)\rangle.$$
(46)

$$P_2 = \langle(c, 0.7, 0.5), (d, 0.5, 0.8), (e, 0.6, 0.6)\rangle, Q_2 = \langle(cd, 0.5, 0.8), (de, 0.5, 0.7)\rangle. \tag{47}$$

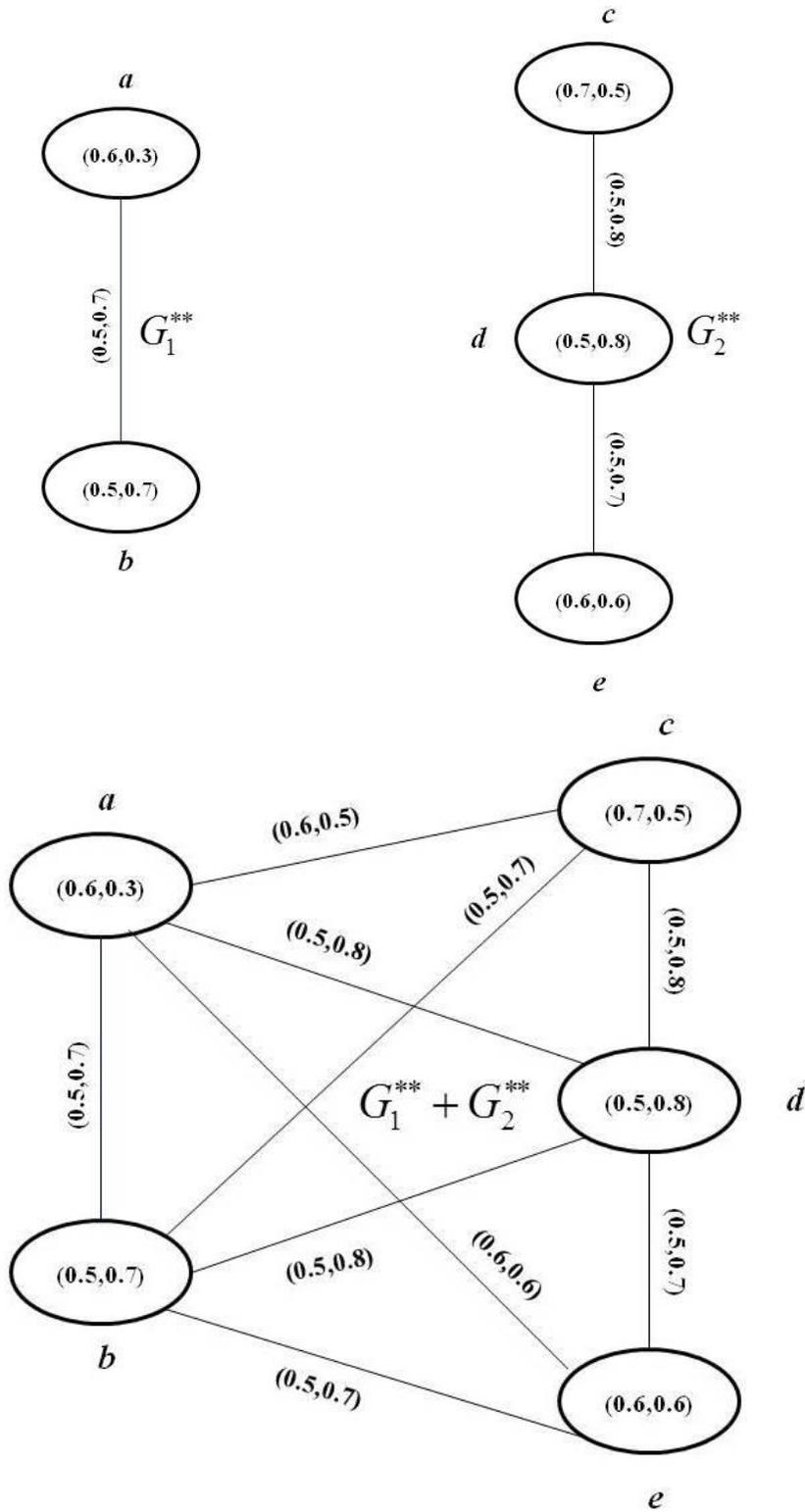

Here $G_1^{**} + G_2^{**}$ represents a Pythagorean fuzzy graph of $G_1 + G_2$.

**Proposition 3.5:** Let $G_1 = (V_1, E_1)$ and $G_2 = (V_2, E_2)$ be crisp graphs and assume $V_1 \cap V_2 \neq \phi$. Further, let $P_1, P_2, Q_1$ and $Q_2$ be Pythagorean fuzzy sets defined on $V_1, V_2, E_1$ and $E_2$, respectively. Then $G_1^{**} \cup G_2^{**} = (P_1 \cup P_2, Q_1 \cup Q_2)$ represents a Pythagorean fuzzy graph of $G = G_1 \cup G_2$ if and only if $G_1^{**} = (P_1, Q_1)$ and $G_2^{**} = (P_2, Q_2)$ are Pythagorean fuzzy graphs of $G_1 = (V_1, E_1)$ and $G_2 = (V_2, E_2)$, respectively.

**Proof:** Suppose that $G_1^{**} \cup G_2^{**} = (P_1 \cup P_2, Q_1 \cup Q_2)$ is a Pythagorean fuzzy graph of $G_1 \cup G_2$. Let $uv \in E_1$, then $uv \notin E_2$ and $u, v \in V_1 - V_2$. Thus

$$\varsigma_{Q_1}(uv) = \varsigma_{Q_1 \cup Q_2}(uv) \leq \min\left(\varsigma_{P_1 \cup P_2}(u), \varsigma_{P_1 \cup P_2}(v)\right)$$
$$= \min\left(\varsigma_{P_1}(u), \varsigma_{P_1}(v)\right). \tag{48}$$

$$\xi_{Q_1}(uv) = \xi_{Q_1 \cup Q_2}(uv) \leq \max\left(\xi_{P_1 \cup P_2}(u), \xi_{P_1 \cup P_2}(v)\right)$$
$$= \max\left(\xi_{P_1}(u), \xi_{P_1}(v)\right). \tag{49}$$

This shows that $G_1^{**} = (P_1, Q_1)$ is a Pythagorean fuzzy graph. Similarly, we can easily prove that $G_1^{**} = (P_2, Q_2)$ is also a Pythagorean fuzzy graph.

This completes the proof. □

**Proposition 3.6:** Let $G_1 = (V_1, E_1)$ and $G_2 = (V_2, E_2)$ be crisp graphs and assume $V_1 \cap V_2 \neq \phi$. Further, let $P_1, P_2, Q_1$ and $Q_2$ be Pythagorean fuzzy sets defined on $V_1, V_2, E_1$ and $E_2$, respectively. Then $G_1^{**} + G_2^{**} = (P_1 + P_2, Q_1 + Q_2)$ is a Pythagorean fuzzy graph of $G = G_1 + G_2$ if and only if $G_1^{**} = (P_1, Q_1)$ and $G_2^{**} = (P_2, Q_2)$ are Pythagorean fuzzy graphs of $G_1 = (V_1, E_1)$ and $G_2 = (V_2, E_2)$, respectively.

**Proof:** The desired result follows from the proof of Propositions 3.4 and 3.5.

**Definition 3.9 *(Complement of Pythagorean Fuzzy Graph)*:** The complement of a Pythagorean fuzzy graph $G^{**} = (P, Q)$ is a Pythagorean fuzzy graph, denoted by $\bar{G}^{**} = (\bar{P}, \bar{Q})$ and defined as follows

(i) $\bar{V} = V$,

(ii) $\begin{cases} \varsigma_{\bar{P}}(u) = \varsigma_P(u), \\ \xi_{\bar{P}}(u) = \xi_P(u). \end{cases} \quad \forall u \in V,$

(iii) $\varsigma_{\bar{Q}}(uv) = \begin{cases} \min(\varsigma_P(u), \varsigma_P(v)), & \text{if } \varsigma_Q(uv) = 0 \\ \min(\varsigma_P(u), \varsigma_P(v)) - \varsigma_Q(uv), & \text{if } \varsigma_Q(uv) > 0 \end{cases}, \quad \forall uv \in E,$

(iv) $\xi_{\bar{Q}}(xy) = \begin{cases} \max(\xi_P(u), \xi_P(v)), & \text{if } \xi_Q(uv) = 0 \\ \max(\xi_P(u), \xi_P(v)) - \xi_P(uv), & \text{if } \xi_Q(uv) > 0 \end{cases}, \quad \forall uv \in E.$

***Example* 3.6:** Let $G^{**} = (P,Q)$ be a Pythagorean fuzzy where

$$P = \langle (a,0.7,0.5),(b,0.3,0.6),(c,0.8,0.2),(d,0.5,0.4) \rangle, \\ Q = \langle (bc,0.3,0.6),(cd,0.4,0.4),(da,0.5,0.4) \rangle \quad . \tag{50}$$

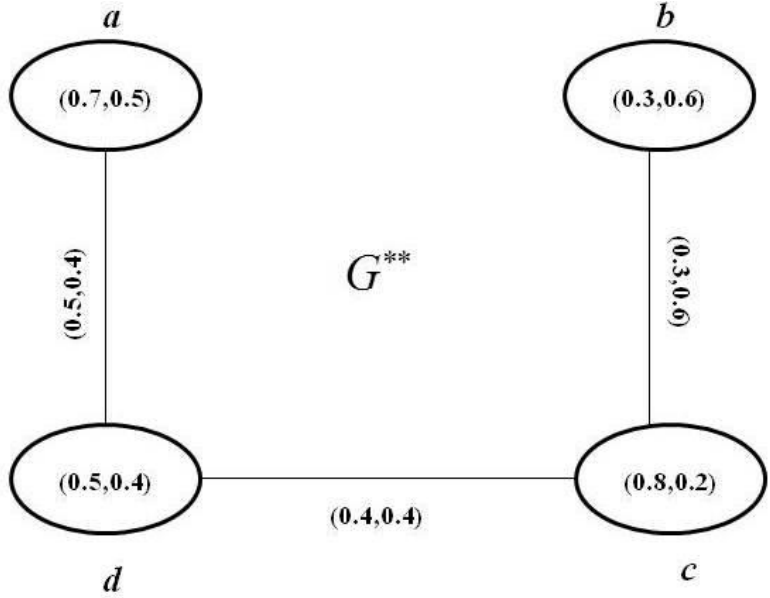

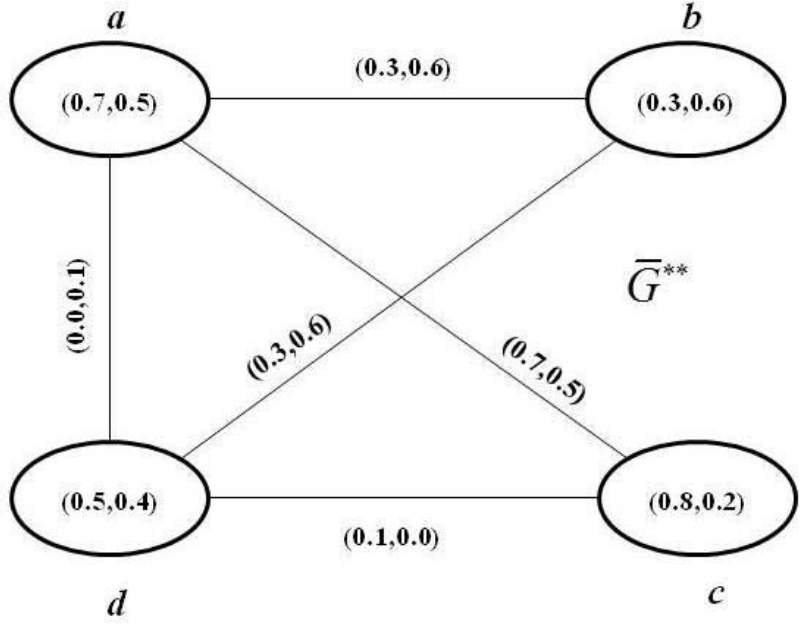

**Note:** One can easily verify that $\bar{\bar{G}}^{**} = G^{**}$.

## 4. Isomorphisms of Pythagorean Fuzzy Graphs

In this section, we characterize various types of (weak) isomorphisms of Pythagorean fuzzy graphs.

***Definition 4.1 (Homomorphism of Pythagorean Fuzzy Graphs):*** Let $G_1^{**} = (P_1, Q_1)$ and $G_2^{**} = (P_2, Q_2)$ be Pythagorean fuzzy graphs. A homomorphism $g$ from $G_1^{**}$ to $G_2^{**}$ is a mapping $g : V_1 \to V_2$ which satisfies the following conditions:

i. $\varsigma_{P_1}(u_1) \leq \varsigma_{P_2}(g(u_1)), \xi_{Q_1}(u_1) \geq \xi_{Q_2}(g(u_1)) \; \forall u_1 \in V_1$,

ii. $\varsigma_{Q_1}(u_1 v_1) \leq \varsigma_{Q_2}(g(u_1) f(v_1)), \xi_{Q_1}(u_1 v_1) \geq \xi_{Q_2}(g(u_1) f(v_1)) \; \forall \; u_1 v_1 \in E_1$.

***Definition 4.2 (Isomorphism of Pythagorean Fuzzy Graphs):*** Let $G_1^{**} = (P_1, Q_1)$ and $G_2^{**} = (P_2, Q_2)$ be Pythagorean fuzzy graphs. An isomorphism $g$ from $G_1^{**}$ to $G_2^{**}$ is a bijective mapping $g : V_1 \to V_2$ which satisfies the following conditions:

i. $\varsigma_{P_1}(u_1) = \varsigma_{P_2}(g(u_1)), \xi_{Q_1}(u_1) = \xi_{Q_2}(g(u_1)) \; \forall u_1 \in V_1$,

ii. $\varsigma_{Q_1}(u_1 v_1) = \varsigma_{Q_2}(g(u_1) f(v_1)), \xi_{Q_1}(u_1 v_1) = \xi_{Q_2}(g(u_1) f(v_1)) \; \forall \; u_1 v_1 \in E_1$.

**Example 3.7:** Let $G_1 = (V_1, E_1)$ and $G_2 = (V_2, E_2)$ be graphs such that $V_1 = \{a_1, a_2, a_3, a_4\}$, $V_2 = \{b_1, b_2, b_3, b_4\}$, $E_1 = \{a_1a_3, a_3a_2, a_4a_1, a_4a_2, a_3a_4\}$ and $E_2 = \{b_1b_2, b_2b_3, b_3b_4, b_4b_1, b_4b_2\}$. Consider two Pythagorean fuzzy graphs $G_1^{**} = (P_1, Q_1)$ and $G_2^{**} = (P_2, Q_2)$, where

$$\left.\begin{array}{l} P_1 = \langle (a_1, 0.4, 0.7), (a_2, 0.7, 0.3), (a_3, 0.6, 0.5), (a_4, 0.3, 0.8) \rangle, \\ Q_1 = \langle (a_1a_3, 0.3, 0.7), (a_3a_2, 0.5, 0.4), (a_2a_4, 0.3, 0.8), (a_4a_1, 0.3, 0.7), (a_3a_4, 0.2, 0.7) \rangle \end{array}\right\}. \quad (51)$$

$$\left.\begin{array}{l} P_2 = \langle (b_1, 0.7, 0.3), (b_2, 0.6, 0.5), (b_3, 0.4, 0.7), (b_4, 0.3, 0.8) \rangle, \\ Q_2 = \langle (b_1b_2, 0.5, 0.4), (b_2b_3, 0.3, 0.7), (b_3b_4, 0.3, 0.7), (b_4b_1, 0.3, 0.8), (b_4b_2, 0.2, 0.7) \rangle \end{array}\right\}. \quad (52)$$

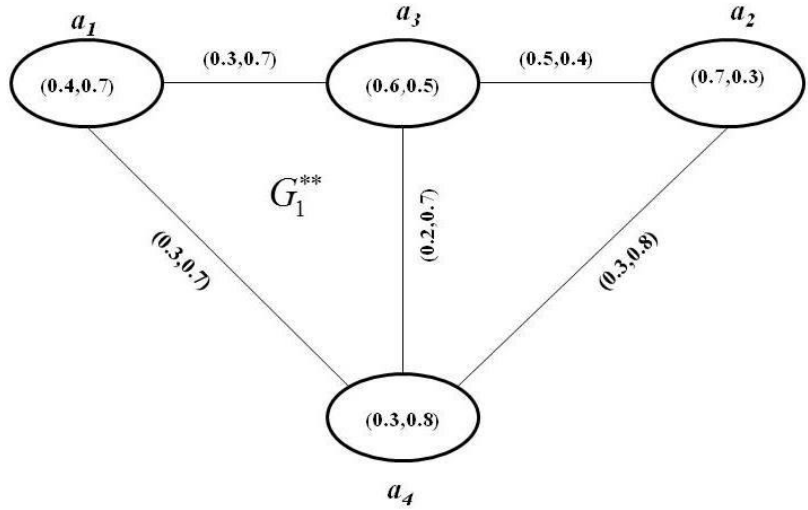

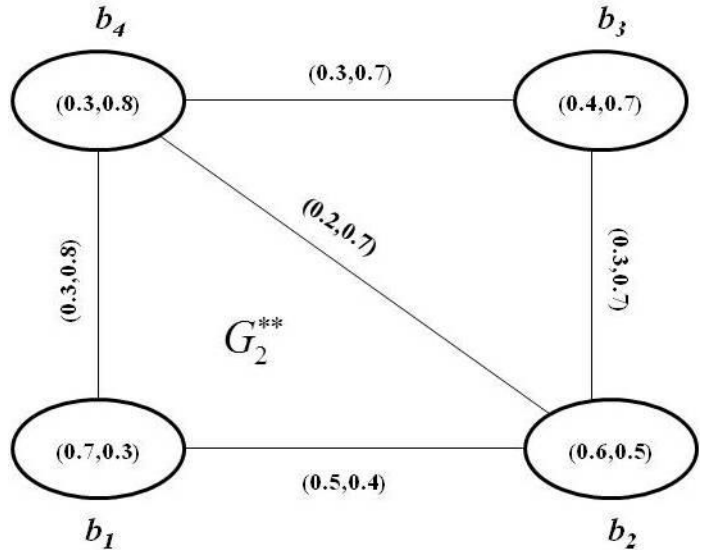

A map $g : V_1 \to V_2$ defined by $g(a_1) = b_3, g(a_2) = b_1, g(a_3) = b_2, g(a_4) = b_4$. Then we have

- $\varsigma_{P_1}(a_1) = \varsigma_{P_2}(g(a_1)), \xi_{P_1}(a_1) = \xi_{P_2}(g(a_1)), \varsigma_{P_1}(a_2) = \varsigma_{P_2}(g(a_2)), \xi_{P_1}(a_2) = \xi_{P_2}(g(a_2)),$
  $\varsigma_{P_1}(a_3) = \varsigma_{P_2}(g(a_3)), \xi_{P_1}(a_3) = \xi_{P_2}(g(a_3)), \varsigma_{P_1}(a_4) = \varsigma_{P_2}(g(a_4)), \xi_{P_1}(a_4) = \xi_{P_2}(g(a_4)).$

- $\varsigma_{Q_1}(a_1 a_3) = \varsigma_{Q_2}(g(a_1)g(a_3)), \xi_{Q_1}(a_1 a_3) = \xi_{Q_2}(g(a_1)g(a_3)), \varsigma_{Q_1}(a_3 a_2) = \varsigma_{Q_2}(g(a_3)g(a_2)),$
  $\xi_{Q_1}(a_3 a_2) = \xi_{Q_2}(g(a_3)g(a_2)), \varsigma_{Q_1}(a_2 a_4) = \varsigma_{Q_2}(g(a_2)g(a_4)), \xi_{Q_1}(a_2 a_4) = \xi_{Q_2}(g(a_2)g(a_4)),$
  $\varsigma_{Q_1}(a_4 a_1) = \varsigma_{Q_2}(g(a_4)g(a_1)), \xi_{Q_1}(a_4 a_1) = \xi_{Q_2}(g(a_4)g(a_1)), \varsigma_{Q_1}(a_3 a_4) = \varsigma_{Q_2}(g(a_3)g(a_4)),$
  $\xi_{Q_1}(a_3 a_4) = \xi_{Q_2}(g(a_3)g(a_4)).$

Hence $G_1^{**}$ is isomorphic to $G_2^{**}$.

***Definition 4.3 (Weak Isomorphism of Pythagorean Fuzzy Graphs):*** Let $G_1^{**} = (P_1, Q_1)$ and $G_2^{**} = (P_2, Q_2)$ be Pythagorean fuzzy graphs. A weak isomorphism $g$ from $G_1^{**}$ to $G_2^{**}$ is a bijective mapping $g: V_1 \to V_2$ which holds the following properties:

i. $g$ is homomorphism.

ii. $\varsigma_{P_1}(u_1) = \varsigma_{P_2}(g(u_1)), \xi_{P_1}(u_1) = \xi_{P_2}(g(u_1)) \ \forall u_1 \in V_1.$

***Example 3.8:*** Let $G_1 = (V_1, E_1)$ and $G_2 = (V_2, E_2)$ be graphs such that $V_1 = \{a_1, a_2\}, V_2 = \{b_1, b_2\}$, $E_1 = \{a_1 a_2\}$ and $E_2 = \{b_1 b_2\}$. Consider two Pythagorean fuzzy graphs $G_1^{**} = (P_1, Q_1)$ and $G_2^{**} = (P_2, Q_2)$, where

$$P_1 = \langle (a_1, 0.5, 0.7), (a_2, 0.6, 0.4) \rangle, Q_1 = \langle (a_1 a_2, 0.5, 0.7) \rangle . \tag{53}$$

$$P_2 = \langle (b_1, 0.6, 0.4), (b_2, 0.5, 0.7) \rangle, Q_2 = \langle (b_1 b_2, 0.4, 0.6) \rangle . \tag{54}$$

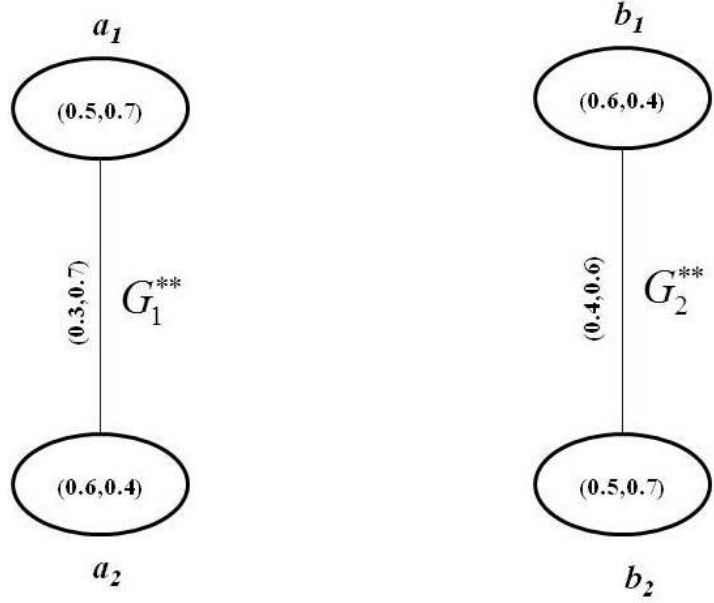

Let $g: V_1 \rightarrow V_2$ be a mapping defined by $g(a_1) = b_2, g(a_2) = b_1$. Then we have

- $\varsigma_{P_1}(a_1) = \varsigma_{P_2}(g(a_1)), \xi_{P_1}(a_1) = \xi_{P_2}(g(a_1)), \varsigma_{P_1}(a_2) = \varsigma_{P_2}(g(a_2)), \xi_{P_1}(a_2) = \xi_{P_2}(g(a_2))$.
- $\varsigma_{Q_1}(a_1 a_2) \neq \varsigma_{Q_2}(g(a_1)f(a_2)), \xi_{Q_1}(a_1 a_2) \neq \xi_{Q_2}(g(a_1)f(a_2))$.

Hence the above mapping represents a weak isomorphism.

***Definition 4.4 (Co-weak Isomorphism of Pythagorean Fuzzy Graphs*** Let $G_1^{**} = (P_1, Q_1)$ and $G_2^{**} = (P_2, Q_2)$ be Pythagorean fuzzy graphs. A co-weak isomorphism $g$ from $G_1^{**}$ to $G_2^{**}$ is a bijective mapping $g: V_1 \rightarrow V_2$ which holds the following properties:

i. $g$ is homomorphism.

ii. $\varsigma_{Q_1}(u_1 v_1) = \varsigma_{Q_2}(g(u_1)f(v_1)), \xi_{Q_1}(u_1 v_1) = \xi_{Q_2}(g(u_1)f(v_1)) \ \forall \ u_1 v_1 \in E_1$.

***Example 3.8:*** Let $G_1 = (V_1, E_1)$ and $G_2 = (V_2, E_2)$ be graphs such that $V_1 = \{a_1, a_2\}, V_2 = \{b_1, b_2\}$, $E_1 = \{a_1 a_2\}$ and $E_2 = \{b_1 b_2\}$. Consider two Pythagorean fuzzy graphs $G_1^{**} = (P_1, Q_1)$ and $G_2^{**} = (P_2, Q_2)$, where

$$P_1 = \langle(a_1, 0.6, 0.5), (a_2, 0.7, 0.6)\rangle, Q_1 = \langle(a_1 a_2, 0.6, 0.5)\rangle . \tag{55}$$

$$P_2 = \langle(b_1, 0.7, 0.6), (b_2, 0.8, 0.5)\rangle, Q_2 = \langle(b_1 b_2, 0.6, 0.5)\rangle . \tag{56}$$

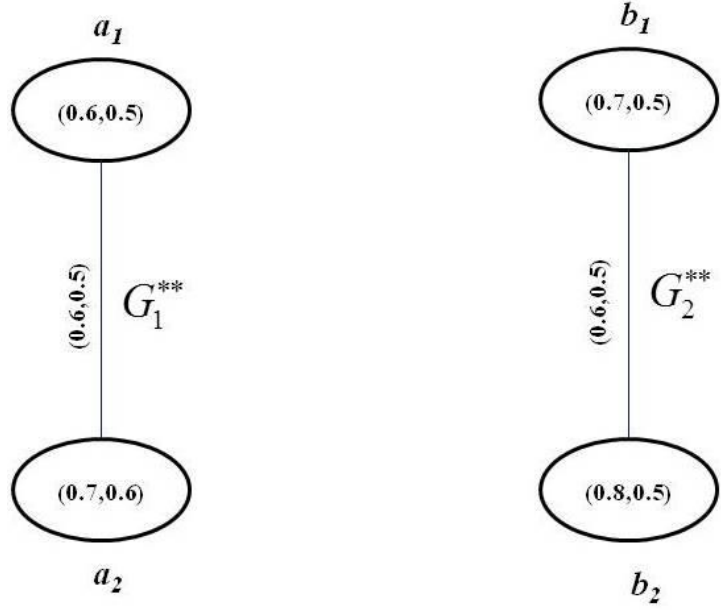

Let $g : V_1 \to V_2$ be a mapping defined by $g(a_1) = b_1, g(a_2) = b_2$. Then we get

- $\varsigma_{P_1}(a_1) \neq \varsigma_{P_2}(g(a_1)), \xi_{P_1}(a_1) \neq \xi_{P_2}(g(a_1)), \varsigma_{P_1}(a_2) \neq \varsigma_{P_2}(g(a_2)), \xi_{P_1}(a_2) \neq \xi_{P_2}(g(a_2))$.
- $\varsigma_{Q_1}(a_1 a_2) = \varsigma_{Q_2}(g(a_1) g(a_2)), \xi_{Q_1}(a_1 a_2) = \xi_{Q_2}(g(a_1) g(a_2))$.

Hence the above mapping represents a co-weak isomorphism.

***Proposition* 4.1:** An isomorphism between Pythagorean fuzzy graphs is an equivalence relation.

***Proof:*** Let $G_1^{**}, G_2^{**}$ and $G_3^{**}$ are three Pythagorean fuzzy graphs. For equivalence relation, we will prove the reflexivity, symmetry, and transitivity for Pythagorean fuzzy graphs.

**i.** *Reflexivity*: It is obvious.

**ii.** *Symmetry*: Let $g : V_1 \to V_2$ be an isomorphism of $G_1^{**}$ onto $G_2^{**}$. Then $g$ is a bijective mapping defined by $g(u_1) = u_2 \ \forall u_1 \in V_1$ and satisfying

$$\varsigma_{P_1}(u_1) = \varsigma_{P_2}(g(u_1)), \xi_{P_1}(u_1) = \xi_{P_2}(g(u_1)) \ \forall \ u_1 \in V_1, \tag{57}$$

$$\varsigma_{Q_1}(u_1 v_1) = \varsigma_{Q_2}(g(u_1) g(v_1)), \xi_{Q_1}(u_1 v_1) = \xi_{Q_2}(g(u_1) g(v_1)) \ \forall \ u_1 v_1 \in E_1. \tag{58}$$

Since $g$ is a bijective mapping then it follows that $g^{-1}(u_2) = u_1 \ \forall u_2 \in V_2$. Thus

$$\varsigma_{P_1}(g^{-1}(u_2)) = \varsigma_{P_2}(u_2), \xi_{P_1}(g^{-1}(u_2)) = \xi_{P_2}(u_2) \ \forall \ u_2 \in V_2, \tag{59}$$

$$\varsigma_{Q_1}(g^{-1}(u_2) g^{-1}(v_2)) = \varsigma_{Q_2}(u_2 v_2), \xi_{Q_1}(g^{-1}(u_2) g^{-1}(v_2)) = \xi_{Q_2}(u_2 v_2) \ \forall \ u_2 v_2 \in E_2. \tag{60}$$

Hence a bijective mapping $g^{-1}: V_2 \rightarrow V_1$ is an isomorphism $G_2^{**}$ onto $G_1^{**}$.

***Transitivity:*** Let $g: V_1 \rightarrow V_2$ and $h: V_2 \rightarrow V_3$ be the isomorphisms of $G_1^{**}$ onto $G_2^{**}$ and $G_2^{**}$ onto $G_3^{**}$, respectively. Since a map $g: V_1 \rightarrow V_2$ defined by $g(u_1) = u_2 \ \forall u_1 \in V_1$ is an isomorphism, so we have

$$\varsigma_{P_1}(u_1) = \varsigma_{P_2}(g(u_1)) = \varsigma_{P_2}(u_2), \xi_{P_1}(u_1) = \xi_{P_2}(g(u_1)) = \xi_{P_2}(u_2) \quad \forall \ u_1 \in V_1, \quad (61)$$

$$\left.\begin{array}{l} \varsigma_{Q_1}(u_1 v_1) = \varsigma_{Q_2}(g(u_1) g(v_1)) = \varsigma_{Q_2}(u_2 v_2), \\ \xi_{Q_1}(u_1 v_1) = \xi_{Q_2}(g(u_1) g(v_1)) = \xi_{Q_2}(v_2 v_2) \quad \forall u_1 v_1 \in E_1 \end{array}\right\}. \quad (62)$$

Similarly, a map $h: V_2 \rightarrow V_3$ defined by $h(u_2) = u_3 \ \forall u_2 \in V_2$ is an isomorphism, so

$$\varsigma_{P_2}(u_2) = \varsigma_{P_3}(h(u_2)) = \varsigma_{P_3}(u_3), \xi_{P_2}(u_2) = \xi_{P_3}(h(u_2)) = \xi_{P_3}(u_3) \quad \forall \ u_2 \in V_2, \quad (63)$$

$$\left.\begin{array}{l} \varsigma_{Q_2}(u_2 v_2) = \varsigma_{Q_3}(h(u_2) h(v_2)) = \varsigma_{Q_3}(u_2 v_2), \\ \xi_{Q_2}(u_2 v_2) = \xi_{Q_3}(h(u_2) h(v_2)) = \xi_{Q_3}(u_2 v_2) \quad \forall \ u_2 v_2 \in E_2 \end{array}\right\}. \quad (64)$$

From (61), (63) and $g(u_1) = u_2 \ \forall u_1 \in V_1$, we have

$$\left.\begin{array}{l} \varsigma_{P_1}(u_1) = \varsigma_{P_2}(g(u_1)) = \varsigma_{P_2}(u_2) = \varsigma_{P_3}(h(u_2)) = \varsigma_{P_3}(h(g(u_1))), \\ \xi_{P_1}(u_1) = \xi_{P_2}(g(u_1)) = \xi_{P_2}(u_2) = \xi_{P_3}(h(u_2)) = \xi_{P_3}(h(g(u_1))) \quad \forall u_1 \in V_1 \end{array}\right\}. \quad (65)$$

From (62) and (64), we have

$$\left.\begin{array}{l} \varsigma_{Q_1}(u_1 v_1) = \varsigma_{Q_2}(g(u_1) g(v_1)) = \varsigma_{Q_2}(u_2 v_2) \\ \qquad = \varsigma_{Q_3}(h(u_2) h(v_2)) = \varsigma_{Q_3}(h(g(u_1)) h(g(v_1))), \\ \xi_{Q_1}(u_1 v_1) = \xi_{Q_2}(g(u_1) g(v_1)) = \xi_{Q_2}(u_2 v_2) \\ \qquad = \xi_{Q_3}(h(u_2) h(v_2)) = \xi_{Q_3}(h(g(u_1)) h(g(v_1))) \quad \forall \ u_1 v_1 \in E_1 \end{array}\right\}. \quad (66)$$

Hence, $g \circ f : V_1 \rightarrow V_3$ is an isomorphism of $G_1^{**}$ onto $G_3^{**}$.

This proves the result. $\square$

***Proposition* 4.2:** The weak isomorphism between Pythagorean fuzzy graphs is a partial order relation.

***Proof*:** Let $G_1^{**}, G_2^{**}$ and $G_3^{**}$ are three Pythagorean fuzzy graphs. For partial order relation, we will prove the reflexivity, anti-symmetry, and transitivity for Pythagorean fuzzy graphs.

i. ***Reflexivity***: It is obvious.

ii. **Anti-Symmetry**: Let $g: V_1 \to V_2$ be a weak isomorphism of $G_1^{**}$ onto $G_2^{**}$. Then $g$ is a bijective mapping defined by $g(u_1) = u_2 \ \forall u_1 \in V_1$ and satisfying

$$\varsigma_{P_1}(u_1) = \varsigma_{P_2}(g(u_1)), \xi_{P_1}(u_1) = \xi_{P_2}(g(u_1)) \ \forall u_1 \in V_1, \tag{67}$$

$$\left.\begin{array}{l} \varsigma_{Q_1}(u_1 v_1) \leq \varsigma_{Q_2}(g(u_1)g(v_1)), \\ \xi_{Q_1}(u_1 v_1) \geq \xi_{Q_2}(g(u_1)g(v_1)) \ \forall u_1 v_1 \in E_1 \end{array}\right\}. \tag{68}$$

Again, let $h: V_2 \to V_1$ be a weak isomorphism of $G_2^{**}$ onto $G_1^{**}$. Then $h$ is a bijective mapping defined by $h(x_2) = x_1 \ \forall x_2 \in V_2$ and satisfying

$$\varsigma_{P_2}(u_2) = \varsigma_{P_1}(h(u_1)), \xi_{P_2}(u_2) = \xi_{P_1}(h(u_1)) \ \forall u_2 \in V_2, \tag{69}$$

$$\left.\begin{array}{l} \varsigma_{Q_2}(u_2 v_2) \leq \varsigma_{Q_1}(h(u_2)h(v_2)), \\ \xi_{Q_2}(u_2 v_2) \geq \xi_{Q_1}(h(u_2)h(v_2)) \ \forall u_2 v_2 \in E_2 \end{array}\right\}. \tag{70}$$

The inequalities (68) and (70) hold on $V_1$ and $V_2$ only when Pythagorean fuzzy graphs $G_1^{**}$ and $G_2^{**}$ have the same numbers of edges and the corresponding edges have the same weight. Hence $G_1^{**}$ and $G_2^{**}$ are identical.

***Transitivity:*** Let $g: V_1 \to V_2$ and $h: V_2 \to V_3$ be the weak isomorphisms of $G_1^{**}$ onto $G_2^{**}$ and $G_2^{**}$ onto $G_3^{**}$, respectively. Since a map $g: V_1 \to V_2$ defined by $g(u_1) = u_2 \ \forall u_1 \in V_1$ is a weak isomorphism, so we have

$$\varsigma_{P_1}(u_1) = \varsigma_{P_2}(g(u_1)) = \varsigma_{P_2}(u_2), \xi_{P_1}(u_1) = \xi_{P_2}(g(u_1)) = \xi_{P_2}(u_2) \ \forall u_1 \in V_1, \tag{71}$$

$$\left.\begin{array}{l} \varsigma_{Q_1}(u_1 v_1) \leq \varsigma_{Q_2}(g(u_1)f(v_1)) = \varsigma_{Q_2}(u_2 v_2), \\ \xi_{Q_1}(u_1 v_1) \geq \xi_{Q_2}(g(u_1)f(v_1)) = \xi_{Q_2}(v_2 v_2) \ \forall u_1 v_1 \in E_1 \end{array}\right\}. \tag{72}$$

Similarly, a map $h: V_2 \to V_3$ defined by $h(u_2) = u_3 \ \forall u_2 \in V_2$ is a weak isomorphism, so

$$\varsigma_{P_2}(u_2) = \varsigma_{P_3}(h(u_2)) = \varsigma_{P_3}(u_3), \xi_{P_2}(u_2) = \xi_{P_3}(h(u_2)) = \xi_{P_3}(u_3) \ \forall u_2 \in V_2, \tag{73}$$

$$\left.\begin{array}{l} \varsigma_{Q_2}(u_2 v_2) \leq \varsigma_{Q_3}(h(u_2)g(v_2)) = \varsigma_{Q_3}(u_2 v_2), \\ \xi_{Q_2}(u_2 v_2) \geq \xi_{Q_3}(h(u_2)g(v_2)) = \xi_{Q_3}(u_2 v_2) \ \forall u_2 v_2 \in E_2 \end{array}\right\}. \tag{74}$$

From (71), (73) and $g(u_1) = u_2 \quad \forall u_1 \in V_1$, we have

$$\left.\begin{array}{l} \varsigma_{P_1}(u_1) = \varsigma_{P_2}(g(u_1)) = \varsigma_{P_2}(u_2) = \varsigma_{P_3}(h(u_2)) = \varsigma_{P_3}(h(g(u_1))), \\ \xi_{P_1}(u_1) = \xi_{P_2}(g(u_1)) = \xi_{P_2}(u_2) = \xi_{P_3}(h(u_2)) = \xi_{P_3}(h(g(u_1))) \quad \forall u_1 \in V_1 \end{array}\right\}. \quad (75)$$

From (72) and (74), we have

$$\left.\begin{array}{l} \varsigma_{Q_1}(u_1 v_1) \leq \varsigma_{Q_2}(g(u_1) g(v_1)) = \varsigma_{Q_2}(u_2 v_2) \\ \qquad = \varsigma_{Q_3}(h(u_2) h(v_2)) = \varsigma_{Q_3}(h(g(u_1)) h(g(v_1))), \\ \xi_{Q_1}(u_1 v_1) \geq \xi_{Q_2}(g(u_1) g(v_1)) = \xi_{Q_2}(u_2 v_2) \\ \qquad = \xi_{Q_3}(h(u_2) h(v_2)) = \xi_{Q_3}(h(g(u_1)) h(g(v_1))) \quad \forall u_1 v_1 \in E_1 \end{array}\right\}. \quad (76)$$

Hence, $h \circ g : V_1 \to V_3$ is a weak isomorphism of $G_1^{**}$ onto $G_3^{**}$.

This proves the Proposition. □

***Definition* 4.5:** A Pythagorean fuzzy graph $G^{**}(P,Q)$ is called self-complementary if $G^{**} \cong \overline{G}^{**}$.

***Proposition* 4.3:** If $G^{**}(P,Q)$ is a self-complementary Pythagorean fuzzy graph, then

$$\sum_{u \neq v} \varsigma_Q(uv) = \frac{1}{2} \sum_{u \neq v} \min(\varsigma_P(u), \varsigma_P(v)), \quad (77)$$

$$\sum_{u \neq v} \xi_Q(uv) = \frac{1}{2} \sum_{u \neq v} \max(\xi_P(u), \xi_P(v)). \quad (78)$$

***Proof*:** Let $G^{**} = (P,Q)$ be a self-complementary PFG. Then there exist isomorphisms such that

$$\varsigma_{\overline{P}}(g(u)) = \varsigma_P(u) \text{ and } \xi_{\overline{P}}(g(u)) = \xi_P(u) \quad \forall u \in V, \quad (79)$$

and $\quad \xi_{\overline{Q}}(g(u)g(v)) = \xi_Q(uv) \text{ and } \xi_{\overline{Q}}(g(u)g(v)) = \xi_Q(uv) \quad \forall u,v \in V. \quad (80)$

By the definition of $\overline{G}^{**}$, we have

$$\varsigma_{\overline{Q}}(g(u)g(v)) = \min(\varsigma_{\overline{P}}(g(u)), \varsigma_{\overline{P}}(gv)) - \varsigma_Q(g(u)g(v)) \quad (81)$$

i.e. $\qquad \varsigma_Q(uv) = \min(\varsigma_P(u), \varsigma_P(v)) - \varsigma_Q(g(u)g(v))$

i.e. $\qquad \sum_{u \neq v} \varsigma_Q(uv) + \sum_{u \neq v} \varsigma_Q(g(u)g(v)) = \sum_{u \neq v} \min(\varsigma_P(u), \varsigma_P(v))$

i.e. $\qquad 2\sum_{u \neq v} \varsigma_Q(uv) = \sum_{u \neq v} \min(\varsigma_P(u), \varsigma_P(v))$

i.e.
$$\sum_{u \neq v} \varsigma_Q(uv) = \frac{1}{2} \sum_{u \neq v} \min(\varsigma_P(u), \varsigma_P(v)) \tag{82}$$

Similarly, we can show that

$$\sum_{u \neq v} \xi_Q(uv) = \frac{1}{2} \sum_{u \neq v} \min(\xi_P(u), \xi_P(v)) \tag{83}$$

This proves the proposition. □

*Remark*: The conditions given in Proposition 4.3 are not sufficient. In the following example $G^{**}$ is not isomorphic to $\bar{G}^{**}$ but

$$\sum_{u \neq v} \varsigma_Q(uv) = 0.3 + 0.15 + 0.1 = 0.55, \sum_{u \neq v} \xi_Q(uv) = 0.4 + 0.3 + 0.3 = 1.0$$

$$\frac{1}{2} \sum_{u \neq v} \min(\varsigma_P(u), \varsigma_P(v)) = \frac{1}{2}(0.5 + 0.3 + 0.3) = 0.55,$$

$$\frac{1}{2} \sum_{u \neq v} \min(\xi_P(u), \xi_P(v)) = \frac{1}{2}(0.7 + 0.7 + 0.6) = 1.0.$$

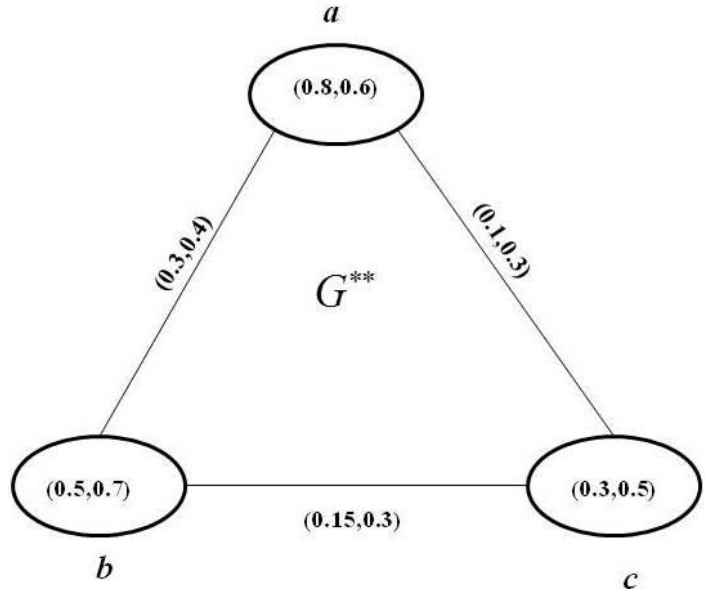

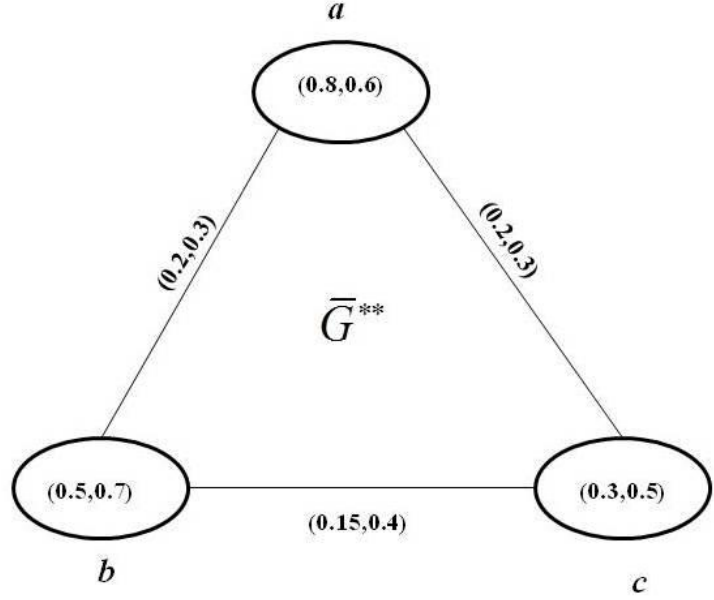

**Proposition 4.4:** For any Pythagorean fuzzy graph $G^{**} = (P, Q)$, if

$$\varsigma_Q(uv) = \frac{1}{2} \min(\varsigma_P(u), \varsigma_P(v)) \tag{84}$$

$$\xi_Q(uv) = \frac{1}{2} \max(\xi_P(u), \xi_P(v)) \quad \forall u, v \in V, \tag{85}$$

then $G^{**}$ is a self-complementary Pythagorean fuzzy graph.

**Proof:** Assume $G^{**} = (P, Q)$ is a Pythagorean fuzzy graph such that

$$\varsigma_Q(uv) = \frac{1}{2} \min(\varsigma_P(u), \varsigma_P(v)) \; \xi_Q(uv) = \frac{1}{2} \max(\xi_P(u), \xi_P(v)) \quad \forall u, v \in V.$$

Then $G^{**} \cong \overline{G}^{**}$ under the identity mapping defined on $V$.

This proves the result. □

**Remark 4.1:** The conditions given in Proposition 4.4 are not necessary. In the following example $G^{**} \cong \overline{G}^{**}$, where the isomorphism $f : V \to V$ is given by $f(a) = b, f(b) = d, f(c) = a, f(d) = c, f(e) = e$, but

$$\varsigma_Q(uv) \neq \frac{1}{2} \min(\varsigma_P(u), \varsigma_P(v)), \; \xi_Q(uv) \neq \frac{1}{2} \max(\xi_P(u), \xi_P(v)) \quad \forall u, v \in V$$

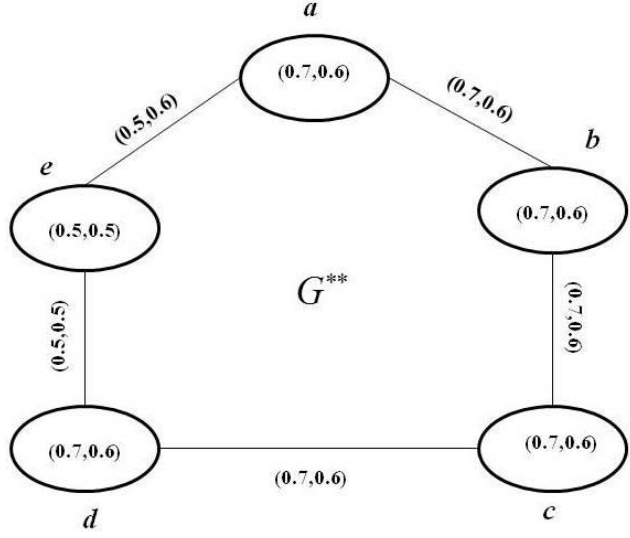

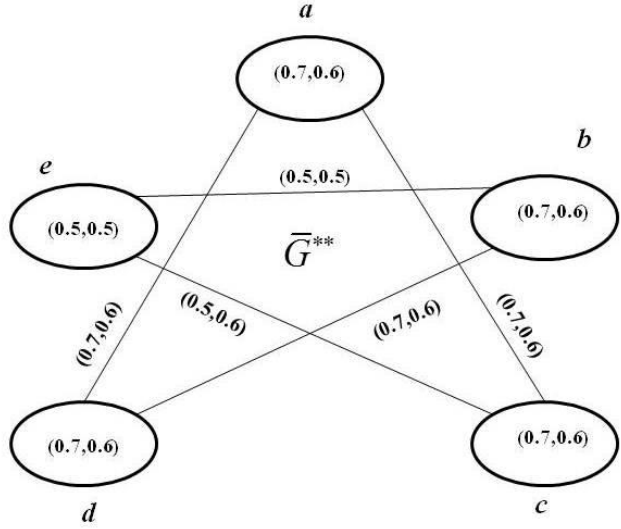

***Proposition* 4.3:** Let $G_1^{**}(P_1, Q_1)$ and $G_2^{**}(P_2, Q_2)$ be two Pythagorean fuzzy graphs. Then

**(i)** $G_1^{**} + G_2^{**} \cong \overline{G}_1^{**} \cup \overline{G}_2^{**}$,

**(ii)** $\overline{G_1^{**} \cup G_2^{**}} \cong \overline{G}_1^{**} + \overline{G}_2^{**}$.

***Proof*:** **(i)** Let $I : V_1 \cup V_2 \to V_1 \cup V_2$ be the identity map. To prove (i), it is enough to prove that

$\varsigma_{\overline{P_1+P_2}}(u) = \varsigma_{\overline{P}_1 \cup \overline{P}_2}(u), \xi_{\overline{P_1+P_2}}(u) = \xi_{\overline{P}_1 \cup \overline{P}_2}(u)$ and $\varsigma_{\overline{Q_1+Q_2}}(uv) = \varsigma_{\overline{Q}_1 \cup \overline{Q}_2}(uv), \xi_{\overline{Q_1+Q_2}}(uv) = \xi_{\overline{Q}_1 \cup \overline{Q}_2}(uv)$.

From the definition of the complement of the Pythagorean fuzzy graph, we have

$\varsigma_{\overline{P_1+P_2}}(u) = \varsigma_{P_1+P_2}(u)$

$$= \begin{cases} \varsigma_{P_1}(u), & \text{if } u \in V_1 \\ \varsigma_{P_2}(u) & \text{if } u \in V_2. \end{cases}$$

$$= \begin{cases} \varsigma_{\bar{P}_1}(u), & \text{if } u \in V_1 \\ \varsigma_{\bar{P}_2}(u) & \text{if } u \in V_2. \end{cases}$$

$$= \varsigma_{\bar{P}_1 \cup \bar{P}_2}(u). \tag{86}$$

$$\xi_{\overline{P_1+P_2}}(u) = \xi_{P_1+P_2}(u)$$

$$= \begin{cases} \xi_{P_1}(u), & \text{if } u \in V_1 \\ \xi_{P_2}(u) & \text{if } u \in V_2. \end{cases}$$

$$= \begin{cases} \xi_{\bar{P}_1}(u), & \text{if } u \in V_1 \\ \xi_{\bar{P}_2}(u) & \text{if } u \in V_2. \end{cases}$$

$$= \xi_{\bar{P}_1 \cup \bar{P}_2}(u). \tag{87}$$

$$\varsigma_{\overline{Q_1+Q_2}}(uv) = \min\left(\varsigma_{Q_1+Q_2}(u), \varsigma_{Q_1+Q_2}(v)\right) - \varsigma_{(Q_1+Q_2)}(uv)$$

$$= \begin{cases} \min\left(\varsigma_{P_1 \cup P_2}(u), \varsigma_{P_1 \cup P_2}(v)\right) - \varsigma_{(Q_1 \cup Q_2)}(uv), & \text{if } uv \in E_1 \bigcup E_2 \\ \min\left(\varsigma_{P_1 \cup P_2}(u), \varsigma_{P_1 \cup P_2}(v)\right) - \min\left(\varsigma_{P_1}(u), \varsigma_{P_2}(v)\right) & \text{if } uv \in E' \end{cases}$$

$$= \begin{cases} \min\left(\varsigma_{P_1}(u), \varsigma_{P_1}(v)\right) - \varsigma_{Q_1}(uv), & \text{if } uv \in E_1 \\ \min\left(\varsigma_{P_2}(u), \varsigma_{P_2}(v)\right) - \varsigma_{Q_2}(uv), & \text{if } uv \in E_2 \\ \min\left(\varsigma_{P_1}(u), \varsigma_{P_2}(v)\right) - \min\left(\varsigma_{P_1}(u), \varsigma_{P_2}(v)\right) & \text{if } uv \in E' \end{cases}$$

$$= \begin{cases} \varsigma_{\bar{Q}_1}(uv), & \text{if } uv \in E_1 \\ \varsigma_{\bar{Q}_2}(uv), & \text{if } uv \in E_2 \\ 0 & \text{if } uv \in E' \end{cases}$$

$$= \varsigma_{\bar{Q}_1 \cup \bar{Q}_2}(uv). \tag{88}$$

$$\xi_{\overline{Q_1+Q_2}}(uv) = \max\left(\xi_{P_1+P_2}(u), \xi_{P_1+P_2}(v)\right) - \xi_{(Q_1+Q_2)}(uv)$$

$$= \begin{cases} \max\left(\xi_{P_1 \cup P_2}(u), \xi_{P_1 \cup P_2}(v)\right) - \xi_{(Q_1 \cup Q_2)}(uv), & \text{if } uv \in E_1 \bigcup E_2 \\ \max\left(\xi_{P_1 \cup P_2}(u), \xi_{P_1 \cup P_2}(v)\right) - \max\left(\xi_{P_1}(u), \xi_{P_2}(v)\right) & \text{if } uv \in E' \end{cases}$$

$$= \begin{cases} \max(\xi_{P_1}(u), \xi_{P_1}(v)) - \xi_{Q_1}(uv), & \text{if } uv \in E_1 \\ \max(\xi_{P_2}(u), \xi_{P_2}(v)) - \xi_{Q_2}(uv), & \text{if } uv \in E_2 \\ \max(\xi_{P_1}(u), \xi_{P_2}(v)) - \max(\xi_{P_1}(u), \xi_{P_2}(v)) & \text{if } uv \in E' \end{cases}$$

$$= \begin{cases} \xi_{\bar{Q}_1}(uv), & \text{if } uv \in E_1 \\ \xi_{\bar{Q}_2}(uv), & \text{if } uv \in E_2 \\ 0 & \text{if } uv \in E' \end{cases}$$

$$= \xi_{\bar{Q}_1 \cup \bar{Q}_2}(uv). \tag{89}$$

This proves the result. □

**(ii)** Let $I: V_1 \cup V_2 \to V_1 \cup V_2$ be the identity map. To prove (ii), it is enough to prove that

$$\varsigma_{\overline{P_1 \cup P_2}}(u) = \varsigma_{\bar{P}_1 + \bar{P}_2}(u), \xi_{\overline{P_1 \cup P_2}}(u) = \xi_{\bar{P}_1 + \bar{P}_2}(u) \text{ and } \varsigma_{\overline{Q_1 \cup Q_2}}(uv) = \varsigma_{\bar{Q}_1 + \bar{Q}_2}(uv), \xi_{\overline{Q_1 \cup Q_2}}(uv) = \xi_{\bar{Q}_1 + \bar{Q}_2}(uv).$$

From the definition of the complement of the Pythagorean fuzzy graph, we have

$$\varsigma_{\overline{P_1 \cup P_2}}(u) = \varsigma_{P_1 \cup P_2}(u)$$

$$= \begin{cases} \varsigma_{P_1}(u), & \text{if } u \in V_1 \\ \varsigma_{P_2}(u) & \text{if } u \in V_2. \end{cases}$$

$$= \begin{cases} \varsigma_{\bar{P}_1}(u), & \text{if } u \in V_1 \\ \varsigma_{\bar{P}_2}(u) & \text{if } u \in V_2. \end{cases}$$

$$= \varsigma_{\bar{P}_1 \cup \bar{P}_2}(u) = \varsigma_{\bar{P}_1 + \bar{P}_2}(u). \tag{90}$$

$$\xi_{\overline{P_1 \cup P_2}}(u) = \xi_{P_1 \cup P_2}(u)$$

$$= \begin{cases} \xi_{P_1}(u), & \text{if } u \in V_1 \\ \xi_{P_2}(u) & \text{if } u \in V_2. \end{cases}$$

$$= \begin{cases} \xi_{\bar{P}_1}(u), & \text{if } u \in V_1 \\ \xi_{\bar{P}_2}(u) & \text{if } u \in V_2. \end{cases}$$

$$= \xi_{\bar{P}_1 \cup \bar{P}_2}(u) = \xi_{\bar{P}_1 + \bar{P}_2}(u). \tag{91}$$

$$\varsigma_{\overline{Q_1 \cup Q_2}}(uv) = \min\left(\varsigma_{P_1 \cup P_2}(u), \varsigma_{P_1 \cup P_2}(v)\right) - \varsigma_{(Q_1 \cup Q_2)}(uv)$$

$$= \begin{cases} \min\left(\varsigma_{P_1}(u), \varsigma_{P_1}(v)\right) - \varsigma_{Q_1}(uv), & \text{if } uv \in E_1 \\ \min\left(\varsigma_{P_2}(u), \varsigma_{P_2}(v)\right) - \varsigma_{Q_2}(uv), & \text{if } uv \in E_2 \\ \min\left(\varsigma_{P_1}(u), \varsigma_{P_2}(v)\right) - 0 & \text{if } uv \in E' \end{cases}$$

$$= \begin{cases} \varsigma_{\overline{Q_1}}(uv), & \text{if } uv \in E_1 \\ \varsigma_{\overline{Q_2}}(uv), & \text{if } uv \in E_2 \\ \min\left(\varsigma_{P_1}(u), \varsigma_{P_2}(v)\right), & \text{if } uv \in E' \end{cases}$$

$$= \begin{cases} \varsigma_{\overline{Q_1} \cup \overline{Q_2}}(uv), & \text{if } uv \in E_1 \text{ or } E_2 \\ \min\left(\varsigma_{P_1}(u), \varsigma_{P_2}(v)\right) & \text{if } uv \in E' \end{cases}$$

$$= \varsigma_{\overline{Q_1} + \overline{Q_2}}(uv). \tag{92}$$

$$\xi_{\overline{Q_1 \cup Q_2}}(uv) = \max\left(\xi_{P_1 \cup P_2}(u), \xi_{P_1 \cup P_2}(v)\right) - \xi_{(Q_1 \cup Q_2)}(uv)$$

$$= \begin{cases} \max\left(\xi_{P_1}(u), \xi_{P_1}(v)\right) - \xi_{Q_1}(uv), & \text{if } uv \in E_1 \\ \max\left(\xi_{P_2}(u), \xi_{P_2}(v)\right) - \xi_{Q_2}(uv), & \text{if } uv \in E_2 \\ \max\left(\xi_{P_1}(u), \xi_{P_2}(v)\right) - 0 & \text{if } uv \in E' \end{cases}$$

$$= \begin{cases} \xi_{\overline{Q_1} \cup \overline{Q_2}}(uv), & \text{if } uv \in E_1 \text{ or } E_2 \\ \max\left(\xi_{P_1}(u), \xi_{P_2}(v)\right) & \text{if } uv \in E' \end{cases}$$

$$= \xi_{\overline{Q_1} + \overline{Q_2}}(uv). \tag{93}$$

This proves the result. □

## 5. Strong and complete Pythagorean Fuzzy Graphs

### 5.1: *Strong Pythagorean Fuzzy Graphs*

***Definition* 5.1.1 ($\mu$-*Strong Pythagorean Fuzzy Graph*):** A PFG $G^{**} = (P, Q)$ is called $\mu$-strong Pythagorean fuzzy graph if

$$\varsigma_Q(uv) = \min\left(\varsigma_P(u), \varsigma_P(v)\right) \text{ and } \xi_Q(uv) \leq \max\left(\xi_P(u), \xi_P(v)\right), \tag{94}$$

and
$$0 \leq \varsigma_Q^2(uv) + \xi_Q^2(uv) \leq 1 \quad \forall uv \in E. \tag{95}$$

*Example*:

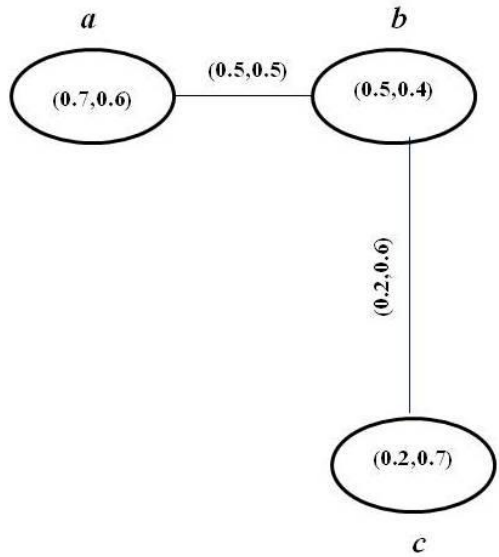

$G^{**}$ is semi μ-strong.

***Definition* 5.1.2 (*v -Strong Pythagorean Fuzzy Graph*):** A PFG $G^{**} = (P,Q)$ is called $v$-strong Pythagorean fuzzy graph if

$$\varsigma_Q(uv) \leq \min(\varsigma_P(u), \varsigma_P(v)) \text{ and } \xi_Q(uv) = \max(\xi_P(u), \xi_P(v)), \quad (96)$$

and
$$0 \leq \varsigma_Q^2(uv) + \xi_Q^2(uv) \leq 1 \quad \forall uv \in E. \quad (97)$$

*Example*:

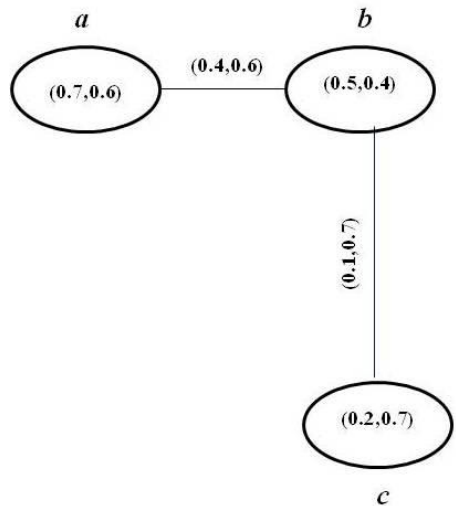

$G^{**}$ is semi v-strong.

***Definition* 5.1.3 (*Strong Pythagorean Fuzzy Graph*):** A PFG $G^{**} = (P, Q)$ is called strong Pythagorean fuzzy graph if

$$\varsigma_Q(uv) = \min(\varsigma_P(u), \varsigma_P(v)) \text{ and } \xi_Q(uv) = \max(\xi_P(u), \xi_P(v)), \tag{98}$$

and
$$0 \leq \varsigma_Q^2(uv) + \xi_Q^2(uv) \leq 1 \quad \forall uv \in E. \tag{99}$$

*Example:*

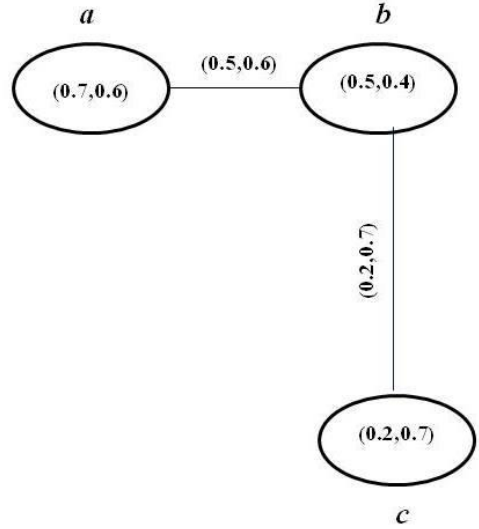

$G^{**}$ is strong.

***Proposition* 5.1.1:** If $G_1^{**}$ and $G_2^{**}$ are strong-PFGs, then $G_1^{**} \times G_2^{**}$, $G_1^{**}[G_2^{**}]$ and $G_1^{**} + G_2^{**}$ are also strong-PFGs.

***Proof*:** These can easily prove on lines similar to proof of Propositions 3.1, 3.2 and 3.4. □

***Remark* 2:** The union of two strong-PFGs is not necessarily a strong-PFG.

***Example*:** We consider the following example:

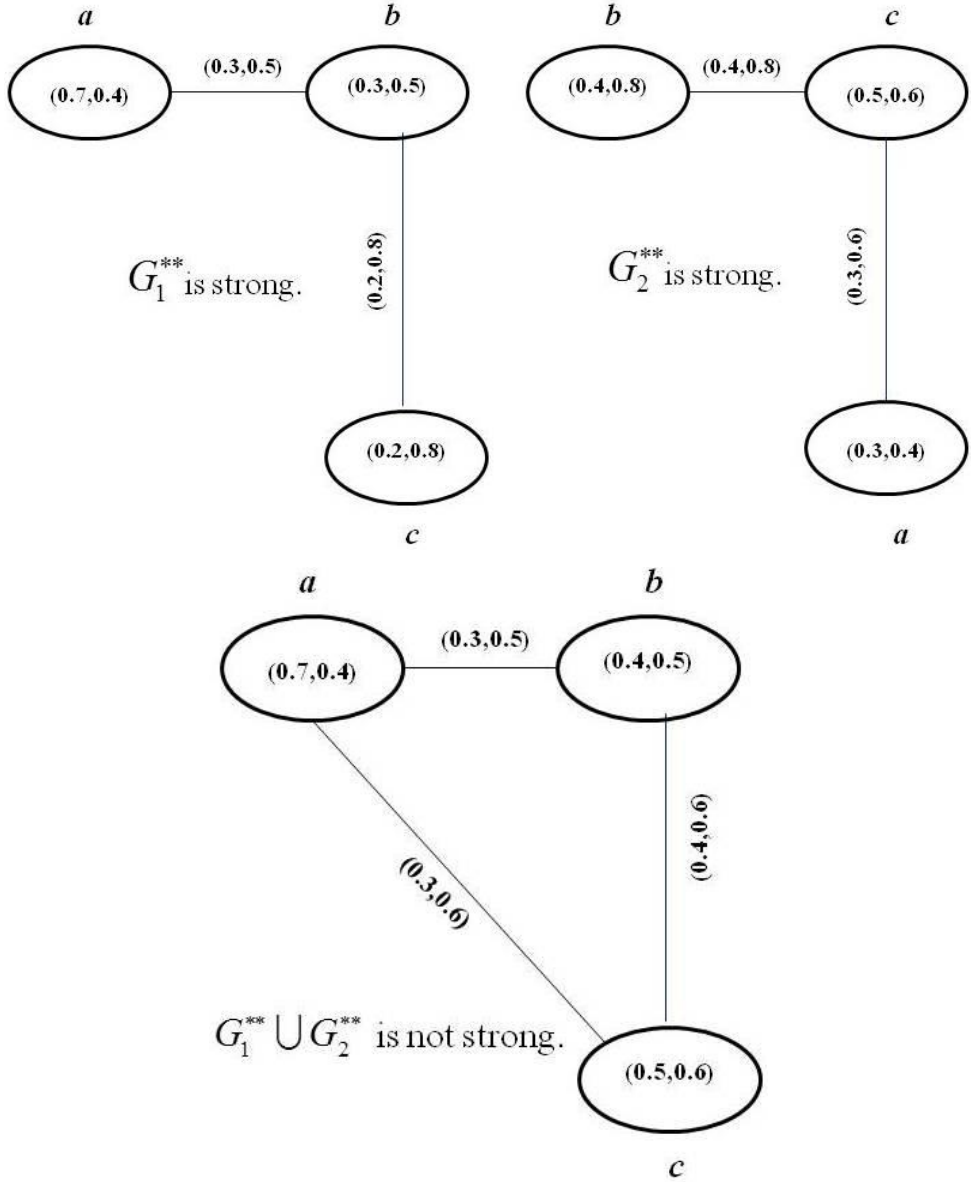

**Proposition 5.1.2:** If $G_1^{**} \times G_2^{**}$ is a strong-PFG, then at least $G_1^{**}$ or $G_2^{**}$ must be strong.

**Proof:** Assume $G_1^{**}$ and $G_2^{**}$ are not strong-PFGs. Then there exists $u_1 v_1 \in E_1$ and $u_2 v_2 \in E_2$ such that

$$\varsigma_{Q_1}(u_1 v_1) < \min(\varsigma_{P_1}(u), \varsigma_{P_1}(v)), \varsigma_{Q_2}(uv) < \min(\varsigma_{P_2}(u), \varsigma_{P_2}(v)), \quad (100)$$

$$\xi_{Q_1}(u_1 v_1) < \max(\xi_{P_1}(u), \xi_{P_1}(v)), \xi_{Q_2}(u_1 v_1) < \max(\xi_{P_2}(u_1), \xi_{P_2}(v_1)). \quad (101)$$

Further, assume that

$$\varsigma_{Q_2}(u_2v_2) \leq \varsigma_{Q_1}(u_1v_1) < \min(\varsigma_{P_1}(u_1), \varsigma_{P_1}(v_1)) \leq \varsigma_{P_1}(u_1). \tag{102}$$

Let $E = \{(t,u_2)(t,v_2) \mid t \in V_1, u_2v_2 \in E_2\} \cup \{(u_1,w)(v_1,w) \mid w \in V_2, u_1v_1 \in E_1\}$ and consider $(t,u_2)(t,v_2) \in E$, then we have

$$\varsigma_{Q_1 \times Q_2}((t,u_2)(t,v_2)) = \min(\varsigma_{P_1}(t), \varsigma_{Q_2}(u_2v_2)) < \min(\varsigma_{P_1}(t), \varsigma_{P_2}(u_2), \varsigma_{P_2}(v_2)), \tag{103}$$

and

$$\varsigma_{P_1 \times P_2}(u_1, u_2) = \min(\varsigma_{P_1}(u_1), \varsigma_{P_2}(u_2)) \text{ and } \varsigma_{P_1 \times P_2}(u_1, v_2) = \min(\varsigma_{P_1}(u_1), \varsigma_{P_2}(v_2)). \tag{104}$$

Therefore,

$$\min(\varsigma_{P_1 \times P_2}(t,u_2), \varsigma_{P_1 \times P_2}(t,v_2)) = \min(\varsigma_{P_1}(t), \varsigma_{P_2}(u_2)\varsigma_{P_2}(v_2)). \tag{105}$$

Hence

$$\varsigma_{Q_1 \times Q_2}((t,u_2)(t,v_2)) < \min(\varsigma_{P_1 \times P_2}(t,u_2), \varsigma_{P_1 \times P_2}(t,v_2)). \tag{106}$$

Similarly, we can easily show that

$$\xi_{Q_1 \times Q_2}((t,u_2)(t,v_2)) < \max(\xi_{P_1 \times P_2}(t,u_2), \xi_{P_1 \times P_2}(t,v_2)). \tag{107}$$

From (106) and (107), $G_1^{**} \times G_2^{**}$ is not a strong-PFG. This is a contradiction.

This proves the Proposition. □

***Remark***: If $G_1^{**}$ is strong and $G_2^{**}$ is not strong, then $G_1^{**} \times G_2^{**}$ may or may not be strong.

***Example* 3.8:** We consider the following examples:

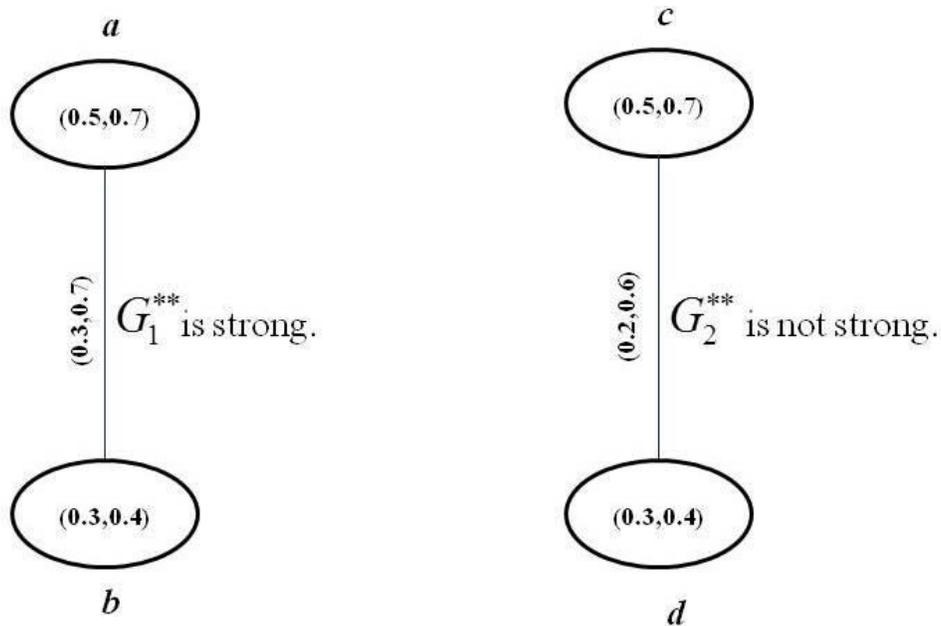

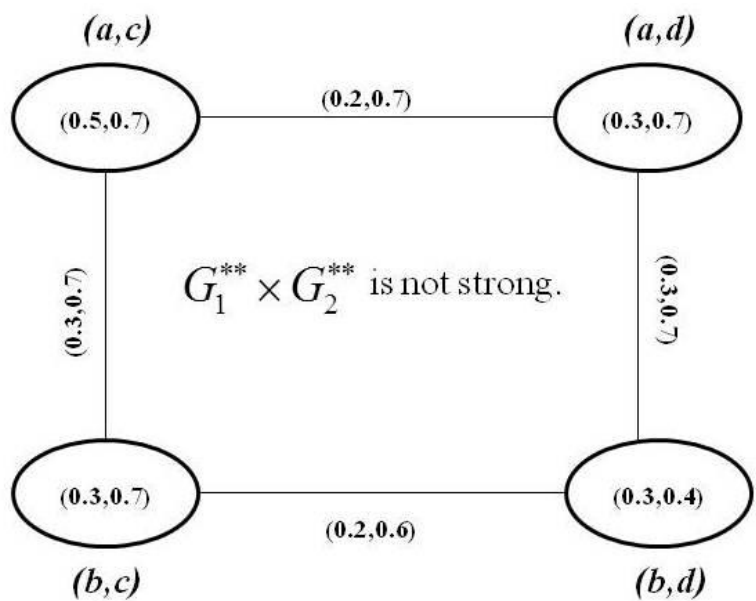

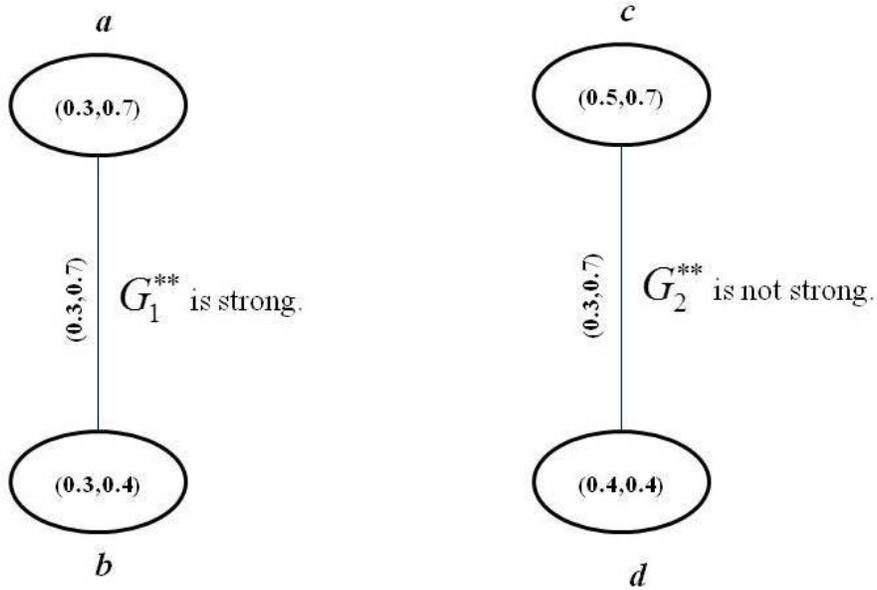

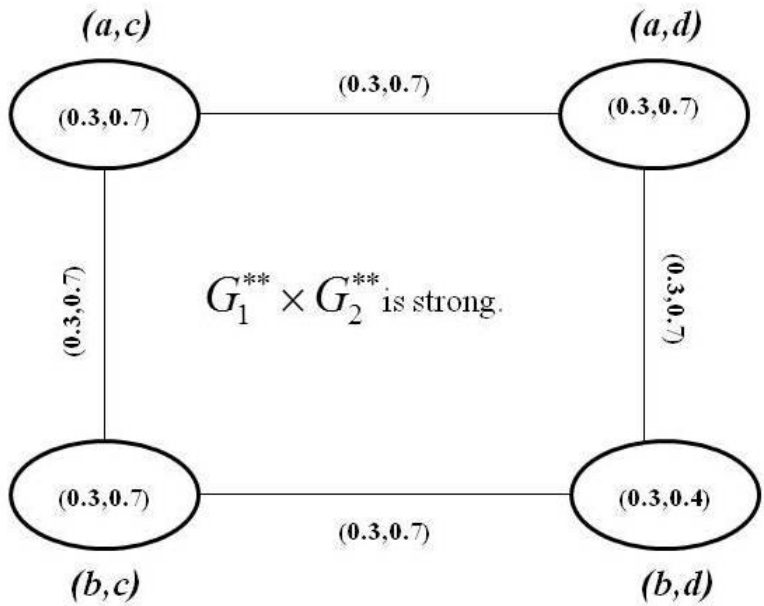

***Proposition* 5.1.3:** If $G_1^{**}\left[G_2^{**}\right]$ is a strong-PFG, then at least $G_1^{**}$ or $G_2^{**}$ must be strong.

***Proof*:** It can be proved in a similar way as Proposition 5.1.2. □

***Proposition* 5.1.4:** If $G_1^{**}+G_2^{**}$ is a strong-PFG, then at least $G_1^{**}$ or $G_2^{**}$ must be strong.

***Proof*:** The proof is similar to Proposition 5.12. □

***Definition* 5.1.4 (*Complement of Strong Pythagorean Fuzzy Graph*):** Let $G^{**} = (P,Q)$ be a strong-PFG of $G = (V,E)$. The complement of $G^{**} = (P,Q)$ is a strong-PFG, denoted by $\overline{G^{**}} = (\overline{P}, \overline{Q})$ and defined as follows

(i) $\overline{V} = V$,

(ii) $\begin{cases} \varsigma_{\overline{P}}(u) = \varsigma_P(u), \\ \xi_{\overline{P}}(u) = \xi_P(u). \end{cases} \quad \forall u \in V$

(iii) $\varsigma_{\overline{Q}}(uv) = \begin{cases} 0, & \text{if } \varsigma_Q(uv) > 0 \\ \min(\varsigma_P(u), \varsigma_P(v)), & \text{if } \varsigma_Q(uv) = 0 \end{cases} \quad \forall uv \in E$

(iv) $\xi_{\overline{Q}}(uv) = \begin{cases} 0, & \text{if } \xi_Q(uv) > 0 \\ \max(\xi_P(u), \xi_P(v)), & \text{if } \xi_Q(uv) = 0 \end{cases} \quad \forall uv \in E$

***Definition* 5.1.5 (*Self Complement of Strong Pythagorean Fuzzy Graph*):** A strong-PFG $G^{**} = (P,Q)$ is called self-complementary if $\overline{\overline{G^{**}}} \cong G^{**}$.

***Example*:** Consider a graph $G = (V,E)$ such as $V = \{a,b,c\}$ and $E = \{ab, bc\}$. Let $G^{**} = (P,Q)$ be a strong-PFG of $G$ where $P = \langle \varsigma_P, \xi_P \rangle$ is a Pythagorean fuzzy set in $V$ and $Q = \langle \varsigma_Q, \xi_Q \rangle$ is a Pythagorean fuzzy set in $E$ defined by

$$P = \langle (a, 0.7, 0.6), (b, 0.5, 0.4), (c, 0.2, 0.3) \rangle, Q = \langle (ab, 0.5, 0.6), (bc, 0.2, 0.7) \rangle . \qquad (108)$$

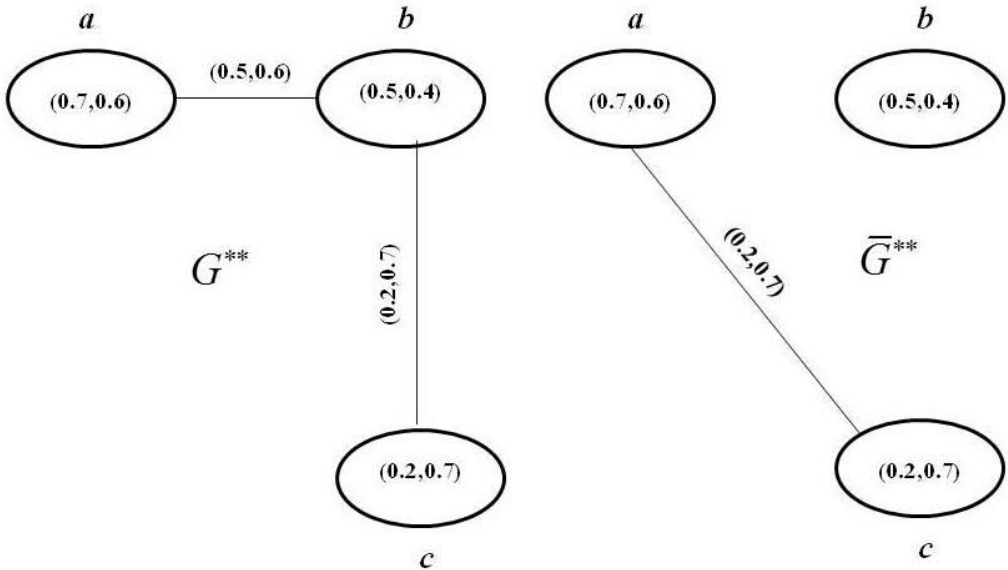

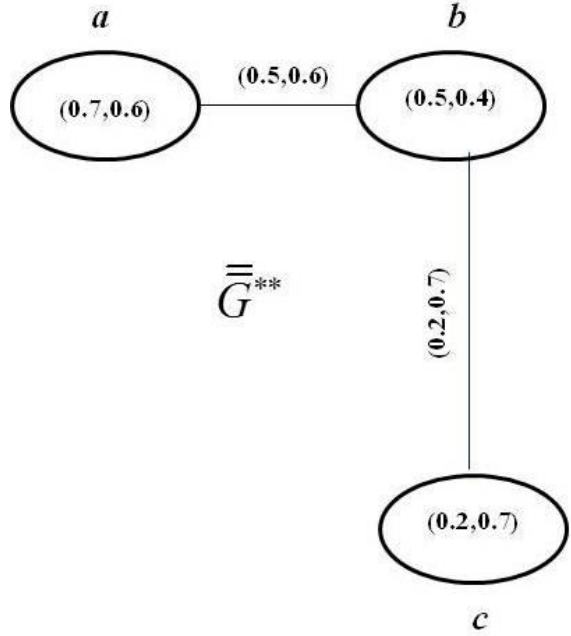

Clearly $\bar{\bar{G}}^{**} \cong G^{**}$, hence $G^{**}$ is a self complementary.

***Proposition 5.1.5:*** Let $G^{**} = (P,Q)$ be a strong-PFG. If $\varsigma_Q(uv) = \min(\varsigma_P(u), \varsigma_P(v))$ and $\xi_Q(uv) = \max(\xi_P(u), \xi_P(v))$ $\forall u, v \in V$, then $G^{**}$ is self-complementary.

***Proof:*** Let $G^{**} = (P,Q)$ be a strong Pythagorean fuzzy graph such that $\varsigma_Q(uv) = \min(\varsigma_P(u), \varsigma_P(v))$ and $\xi_Q(uv) = \max(\xi_P(u), \xi_P(v))$ $\forall u, v \in V$. Identity mapping $I : V \to V$ gives $G^{**} \cong \bar{G}^{**}$. Hence $G^{**}$ is self-complementary. □

***Proposition 5.1.6:*** Let $G^{**} = (P,Q)$ be a self-complementary strong-PFG. Then

$$\sum_{u \neq v} \varsigma_Q(uv) = \sum_{u \neq v} \min(\varsigma_P(u), \varsigma_P(v)), \tag{109}$$

$$\sum_{u \neq v} \xi_Q(uv) = \sum_{u \neq v} \max(\xi_P(u), \xi_P(v)). \tag{110}$$

***Proof:*** Let $G^{**} = (P,Q)$ be a self-complementary strong-PFG. Then there exists an isomorphism $g : V \to V$ such that

$$\varsigma_{\bar{P}}(g(u)) = \varsigma_P(u) \text{ and } \xi_{\bar{P}}(g(u)) = \xi_P(u) \quad \forall u \in V, \tag{111}$$

and $\quad \varsigma_{\bar{Q}}(g(u)g(v)) = \varsigma_Q(uv)$ and $\xi_{\bar{P}}(g(u)g(v)) = \xi_P(uv) \quad \forall u, v \in V.$ (112)

Using the definition of $\bar{G}^{**}$, we have

$$\varsigma_{\bar{Q}}(g(u)g(v)) = \min(\varsigma_{\bar{P}}(g(u)), \varsigma_{\bar{P}}(g(v)))$$

i.e.
$$\varsigma_Q(uv) = \min(\varsigma_P(u), \varsigma_Q(v))$$

or
$$\sum_{u \neq v} \varsigma_Q(uv) = \sum_{u \neq v} \min(\varsigma_P(u), \varsigma_P(v)), \tag{113}$$

Similarly, we can show that

$$\sum_{u \neq v} \xi_Q(uv) = \sum_{u \neq v} \max(\xi_P(u), \xi_P(v)). \tag{114}$$

This completes the proof. □

**Proposition 5.1.7:** Let $G_1^{**}$ and $G_2^{**}$ be two strong-PFGs. Then $G_1^{**} \cong G_2^{**}$ if and only if $\bar{G}_1^{**} \cong \bar{G}_2^{**}$.

**Proof:** Let $G_1^{**}$ and $G_2^{**}$ be isomorphic then there exists a bijective mapping $g: V \to V$ satisfying

$$\varsigma_{P_1}(u) = \varsigma_{P_2}(g(u)) \text{ and } \xi_{P_1}(u) = \xi_{P_2}(g(u)) \quad \forall u \in V_1, \tag{115}$$

and
$$\varsigma_{Q_1}(uv) = \varsigma_{Q_2}(g(u)g(v)) \text{ and } \xi_{Q_1}(uv) = \xi_{Q_2}(g(u)g(v)) \quad \forall uv \in E_1. \tag{116}$$

By the Definition 5.1.4, we have

$$\varsigma_{\bar{Q}_1}(u) = \min(\varsigma_{P_1}(u), \varsigma_{P_1}(y)v) = \min(\varsigma_{P_2}(g(u)), \varsigma_{P_2}(g(v)))$$
$$= \varsigma_{\bar{Q}_2}(g(u)g(v)), \tag{117}$$

and
$$\xi_{\bar{Q}_1}(u) = \min(\xi_{P_1}(u), \xi_{P_1}(y)v) = \max(\xi_{P_2}(g(u)), \xi_{P_2}(g(v)))$$
$$= \xi_{\bar{Q}_2}(g(u)g(v)) \quad \forall uv \in E_1 \tag{118}$$

From (117) and (118), we get $\bar{G}_1^{**} \cong \bar{G}_2^{**}$.

The proof of converse part is straightforward. □

**Proposition 5.1.8:** Let $G_1^{**}$ and $G_2^{**}$ be two strong-PFGs. If there is a weak isomorphism between $G_1^{**}$ and $G_2^{**}$, then there exist a weak isomorphism between $\bar{G}_1^{**}$ and $\bar{G}_2^{**}$.

**Proof:** Let $g$ be a weak isomorphism between $G_1^{**}$ and $G_2^{**}$, then $g: V_1 \to V_2$ is a bijective mapping satisfying $g(u_1) = u_2 \quad \forall u_1 \in V_1$ such that

$$\varsigma_{P_1}(u_1) = \varsigma_{P_2}(g(u_1)) \text{ and } \xi_{P_1}(u_1) = \xi_{P_2}(g(u_1)) \quad \forall u_1 \in V_1, \tag{119}$$

and
$$\varsigma_{Q_1}(u_1v_1) \leq \varsigma_{Q_2}(g(u_1)g(v_1)) \text{ and } \xi_{Q_1}(u_1v_1) \leq \xi_{Q_2}(g(u_1)g(v_1)) \quad \forall u_1v_1 \in E_1. \tag{120}$$

Since $g: V_1 \to V_2$ is a bijective mapping then $g^{-1}: V_2 \to V_1$ is also bijective mapping such that $g^{-1}(u_2) = u_1 \ \forall u_2 \in V_2$. Thus

$$\varsigma_{P_1}(g^{-1}(u_2)) = \varsigma_{P_2}(u_2), \text{ and } \xi_{P_1}(g^{-1}(u_2)) = \xi_{P_2}(u_2) \ \forall u_2 \in V_2. \tag{121}$$

Using the Definition 5.1.4, we have

$$\begin{aligned}
\varsigma_{\bar{Q}_1}(u_1 v_1) &= \min(\varsigma_{P_1}(u_1), \varsigma_{P_1}(v_1)) \\
&\leq \min(\varsigma_{P_2}(g(u_2)), g(v_2))) \\
&= \min(\varsigma_{P_2}(u_2), \varsigma_{P_2}(u_2)) \\
&= \varsigma_{\bar{Q}_2}(u_2 v_2),
\end{aligned} \tag{122}$$

and

$$\begin{aligned}
\xi_{\bar{Q}_1}(u_1 v_1) &= \max(\xi_{P_1}(u_1), \xi_{P_1}(v_1)) \\
&\leq \max(\xi_{P_2}(g(u_2)), g(v_2))) \\
&= \max(\xi_{P_2}(u_2), \xi_{P_2}(u_2)) \\
&= \xi_{\bar{Q}_2}(u_2 v_2),
\end{aligned} \tag{123}$$

Thus $g^{-1}: V_2 \to V_1$ is a weak isomorphism between $\bar{G}_1^{**}$ and $\bar{G}_2^{**}$. □

***Proposition 5.1.9:*** If there is a co-weak isomorphism between two strong-PFGs $G_1^{**}$ and $G_2^{**}$, then there is a homomorphism between $\bar{G}_1^{**}$ and $\bar{G}_2^{**}$.

***Proof:*** It can be proved similar to Proposition 5.1.8. □

**5.2: *Complete Pythagorean Fuzzy Graphs***

***Definition 5.2.1 (Complete $\mu$-Strong Pythagorean Fuzzy Graph):*** A PFG $G^{**} = (P, Q)$ is called a complete $\mu$-strong Pythagorean fuzzy graph if

$$\varsigma_Q(u, v) = \min(\varsigma_P(u), \varsigma_P(v)) \text{ and } \xi_Q(u, v) < \max(\xi_P(u), \xi_P(v)), \tag{124}$$

and
$$0 \leq \varsigma_Q^2(u, v) + \xi_Q^2(u, v) \leq 1 \ \forall \ u, v \in V. \tag{125}$$

***Definition 5.2.2 (Complete $\nu$-Strong Pythagorean Fuzzy Graph):*** A PFG $G^{**} = (P, Q)$ is called a complete $\nu$-strong Pythagorean fuzzy graph if

$$\varsigma_Q(u,v) < \min(\varsigma_P(u), \varsigma_P(v)) \text{ and } \xi_Q(u,v) = \max(\xi_P(u), \xi_P(v)), \quad (126)$$

and
$$0 \leq \varsigma_Q^2(u,v) + \xi_Q^2(u,v) \leq 1 \quad \forall\ u,v \in V. \quad (127)$$

***Definition 5.2.3 (Complete Pythagorean Fuzzy Graph):*** A PFG $G^{**} = (P,Q)$ is called a complete Pythagorean fuzzy graph if

$$\varsigma_Q(u,v) = \min(\varsigma_P(u), \varsigma_P(v)) \text{ and } \xi_Q(u,v) = \max(\xi_P(u), \xi_P(v)), \quad (128)$$

and
$$0 \leq \varsigma_Q^2(u,v) + \xi_Q^2(u,v) \leq 1 \quad \forall\ u,v \in V. \quad (129)$$

**Example:** We consider the following examples:

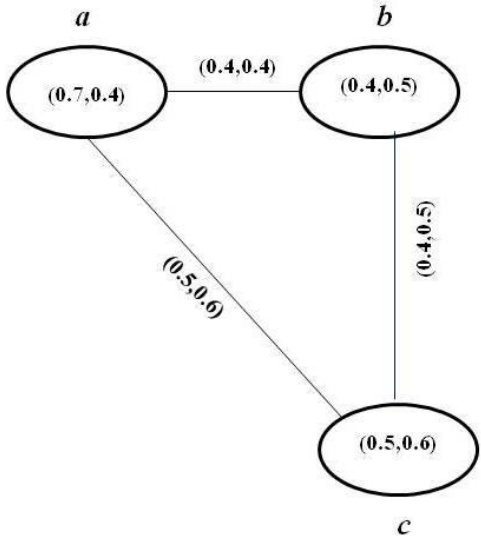
$G^{**}$ is complete μ-strong.

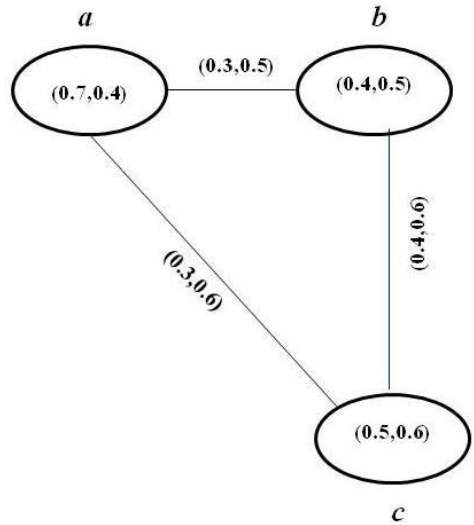
$G^{**}$ is complete ν-strong.

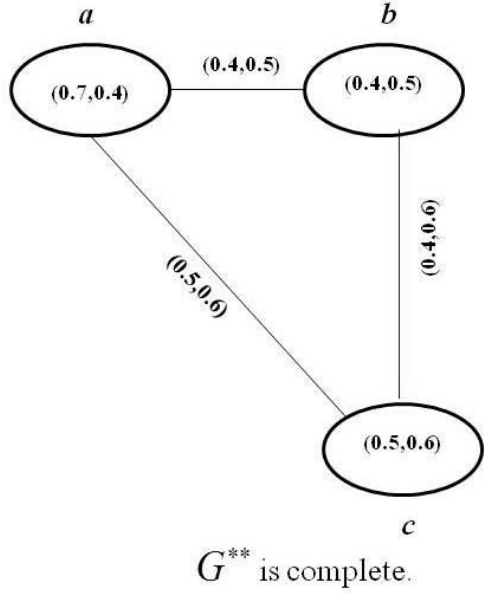

$G^{**}$ is complete.

***Proposition 5.2.1:*** If $G^{**}=(P,Q)$ is a complete-PFG, then $G^{**}[G^{**}]$ is also a complete Pythagorean fuzzy graph.

***Proof:*** It follows from Proposition 3.2. □

***Definition 5.2.4 (Complement of Complete Pythagorean Fuzzy Graph):*** Let $G^{**}=(P,Q)$ be a complete-PFG of $G=(V,E)$. The complement of $G^{**}=(P,Q)$ is a complete Pythagorean fuzzy graph, denoted by $\overline{G}^{**}=(\overline{P},\overline{Q})$ and defined as follows

(v) $\overline{V}=V$,

(vi) $\begin{cases} \varsigma_{\overline{P}}(u)=\varsigma_P(u), \\ \xi_{\overline{P}}(u)=\xi_P(u). \end{cases} \quad \forall x \in V$

(vii) $\varsigma_{\overline{Q}}(uv)=\begin{cases} 0, & \text{if } \varsigma_Q(uv)>0 \\ \min(\varsigma_P(u),\varsigma_P(v)), & \text{if } \varsigma_Q(uv)=0 \end{cases} \quad \forall u,v \in V$

(viii) $\xi_{\overline{Q}}(uv)=\begin{cases} 0, & \text{if } \xi_Q(uv)>0 \\ \max(\xi_P(u),\xi_P(v)), & \text{if } \xi_Q(uv)=0 \end{cases} \quad \forall u,v \in V$

***Definition 5.2.5 (Self Complement of Complete Pythagorean Fuzzy Graph):*** A complete-PFG $G^{**}=(P,Q)$ is called self-complementary if $G^{**} \cong \overline{G}^{**}$.

***Proposition* 5.2.2:** In a self-complementary complete Pythagorean fuzzy graph $G^{**} = (P,Q)$, we have

$$\sum_{u \neq v} \varsigma_Q(uv) = \sum_{u \neq v} \min(\varsigma_P(u), \varsigma_P(v)), \tag{130}$$

and

$$\sum_{u \neq v} \xi_Q(uv) = \sum_{u \neq v} \max(\xi_P(u), \xi_P(v)). \tag{131}$$

**Proof:** The proof is similar to Proposition 5.1.5. □

***Proposition* 5.2.3:** Let $G^{**} = (P,Q)$ be a complete-PFG of $G = (V,E)$. If

$$\varsigma_Q(uv) = \min(\varsigma_P(u), \varsigma_P(v)) \text{ and } \xi_Q(uv) = \max(\xi_P(u), \xi_P(v)) \quad \forall u, v \in V, \tag{132}$$

then $G^{**}$ is self-complementary.

**Proof:** Let $G^{**} = (P,Q)$ be a complete-PFG such that

$$\varsigma_Q(uv) = \min(\varsigma_P(u), \varsigma_P(v)) \text{ and } \xi_Q(uv) = \max(\xi_P(u), \xi_P(v)) \quad \forall u, v \in V,$$

Then $G^{**} = \overline{G}^{**}$ under the identity map, therefore $G^{**}$ is self-complementary. □

***Proposition* 5.2.4:** Let $G_1^{**} = (P_1, Q_1)$ and $G_2^{**} = (P_2, Q_2)$ be complete-PFGs. Then $G_1^{**} \cong G_2^{**}$ if and only if $\overline{G}_1^{**} \cong \overline{G}_2^{**}$.

**Proof:** The proof is similar to Proposition 5.1.7. □

## 6. Conclusions

The work has extended the graph theoretical results under Pythagorean fuzzy environment. We have developed the concept of PFGs as a generalization of fuzzy and intuitionistic fuzzy graphs. Pythagorean fuzzy graphical models provide more precision, flexibility, and compatibility to the user for describing the uncertainty in many combinatorial problems in different areas. We have defined some basic operations such as the Cartesian product, composition, union, join and complement on Pythagorean fuzzy graphs and proved a number of their properties. Further, the work has developed the idea of isomorphism between Pythagorean fuzzy graphs and illustrated with a numerical example. We have also introduced the idea of strong Pythagorean fuzzy graphs and complete Pythagorean fuzzy graphs and proved some results with these graphs.

In future work, we will study different types of graphs including line graphs, hypergraphs, cographs, constant graphs, cycles with Pythagorean fuzzy information. Different types of arcs in Pythagorean fuzzy graphs will also be considered. We will also focus on the applications of Pythagorean fuzzy graphs in different real-life problems.

**Acknowledgments**

Financial support from the Chilean Government (Conicyt) through the Fondecyt Postdoctoral program (Project number-3170556) and Fondecyt Regular program (Project number-1160286) is thankfully acknowledged.

**Conflict of Interests**

The authors declare that there is no conflict of interests regarding the publication of this paper.